\documentclass[11pt,a4paper]{article}
\usepackage[T1]{fontenc}
\usepackage[utf8]{inputenc}
\usepackage[english]{babel}
\usepackage{lmodern,microtype}
\usepackage[margin=26mm,headheight=14pt]{geometry}
\usepackage{amsmath,amssymb,amsthm,mathtools,mathrsfs}
\usepackage{booktabs,longtable,tabularx,array,enumitem}
\usepackage[dvipsnames]{xcolor}
\usepackage{fancyhdr}
\usepackage{hyperref}
\hypersetup{colorlinks=true,linkcolor=MidnightBlue,citecolor=MidnightBlue,urlcolor=MidnightBlue,
	pdftitle={Spectral Fields and Unitary Descent},
	pdfauthor={\'Angel Cano, Nikolay Gusevskii, Marlon M. L\'opez-Flores},
	pdfsubject={Local splitting, simultaneous norms, and transcendental obstructions},
	pdfkeywords={spectral fields, unitary groups, Hermitian forms, norm equations, Zariski density, local fields}}
\newtheorem{theorem}{Theorem}[section]
\newtheorem{lemma}[theorem]{Lemma}
\newtheorem{proposition}[theorem]{Proposition}
\newtheorem{corollary}[theorem]{Corollary}
\theoremstyle{definition}

\theoremstyle{remark}
\newtheorem{remark}[theorem]{Remark}

\numberwithin{equation}{section}
\DeclareMathOperator{\SU}{SU}
\DeclareMathOperator{\SO}{SO}
\DeclareMathOperator{\SL}{SL}
\DeclareMathOperator{\GL}{GL}

\DeclareMathOperator{\tr}{tr}
\DeclareMathOperator{\diag}{diag}

\DeclareMathOperator{\Lie}{Lie}

\DeclareMathOperator{\End}{End}

\newcommand{\Q}{\mathbb Q}
\newcommand{\R}{\mathbb R}
\newcommand{\C}{\mathbb C}

\newcolumntype{Y}{>{\raggedright\arraybackslash}X}
\newcolumntype{L}[1]{>{\raggedright\arraybackslash}p{#1}}
\setlist{itemsep=3pt,topsep=5pt}
\title{\textbf{Spectral Fields and Unitary Descent}\\[3pt]
	\large Local splitting, simultaneous norms, and transcendental obstructions}
\author{\'Angel Cano\qquad Nikolay Gusevskii\qquad Marlon M. L\'opez-Flores}
\date{}
\makeatletter
\renewcommand{\@maketitle}{%
	\begin{center}\vspace*{3mm}{\LARGE\@title\par}
		\vspace{3mm}{\normalsize\@author\par}\end{center}\vspace{3mm}}
\makeatother
\begin{document}
	\maketitle
	\begin{abstract}
		Spectral fields can recover matrix coefficients without standardizing
		the invariant Hermitian form. Over an algebraic subfield of $\C$ stable
		under nontrivial complex conjugation, we prove that every Zariski-dense
		subgroup in a normalized Hermitian model of dimension $m\ge4$ and
		signature $(m-1,1)$ at the fixed embedding contains a loxodromic element
		whose spectrum makes all diagonal coefficients norms. The element can
		be chosen so that every positive power is regular and generates the
		same splitting field. For an arbitrary initial Hermitian model over
		such a field, two spectra suffice for standard unitary descent over
		an extension of relative degree at most $(m!)^2$. Over every quadratic
		extension $F/k$ of nonarchimedean local fields of characteristic zero,
		a four-dimensional quasi-split unitary group of Witt index two
		contains a discrete, Zariski-dense free subgroup whose spectra all
		split over $F$. Over intrinsic fields of positive transcendence degree,
		we construct a single convex cocompact, Zariski-dense free subgroup
		of $\SU(3,1)$ for which, for every initial regular element and every
		normalized extraction, the Hermitian obstruction persists after
		adjoining any finite family of spectra and over the total spectral field.
	\end{abstract}
	
	
	
	\section{Introduction and main results}
	\label{sec:introduction}
	
	Traces and eigenvalues encode arithmetic information about a linear group,
	but recovering its matrix entries over a field and recovering a prescribed
	Hermitian form over the same field are distinct descent problems.
	This distinction is essential for complex hyperbolic groups: a linear
	realization may preserve a Hermitian form that is not similar to the
	standard form over the field of realization. We ask whether adjoining
	spectra of elements of the prescribed group makes its Hermitian form
	similar to the standard model.
	
	The main positive conclusion is standard unitary descent over algebraic
	coefficient fields: two spectra suffice for an arbitrary initial model
	(Corollary~\ref{cor:two-spectrum-algebraic}). Its arithmetic core is
	Theorem~\ref{thm:norm-numeric}, which selects one element making all
	coefficients of a normalized form norms simultaneously.
	The main negative conclusion is a single group over intrinsic
	transcendental fields for which the obstruction survives even the total
	spectral field (Theorem~\ref{thm:transcendental-obstruction} and
	Corollary~\ref{cor:no-total-unitary-descent}). These results determine
	an algebraic range of validity and exhibit a failure beyond it; they
	do not classify all transcendental fields.
	Theorem~\ref{thm:split-local-spectra} is a separate local result that
	tests a proposed route to descent.
	
	This question continues a program on trace fields and fields of definition.
	In complex dimension two, Cunha--Gusevskii~\cite[Theorem~A]{CunhaGusevskii2014}
	realize an irreducible subgroup of $\SU(2,1)$, by unitary conjugation,
	over its trace field with the modulus of an expanding eigenvalue
	adjoined; discreteness is not required. Genzmer~\cite[Theorem~2.3]{Genzmer2014}
	obtains a trace-field realization for a Zariski-dense subgroup of
	$\SU(2,1)$ containing a parabolic element. 
	
	The immediate point of departure is
	Cano--Castro--Gusevskii~\cite[Theorems~1.2--1.3]{CanoCastroGusevskii}:
	traces and one regular loxodromic spectrum give a linear model, while
	standard unitary descent requires the invariant Hermitian form to be
	similar to the standard form. That paper identifies a genuine
	one-spectrum obstruction and proves arithmetic and conditional descent
	results. We study the remaining question: which spectra of elements
	of the prescribed group eliminate its Hermitian norm classes?
	The article can be read independently of those earlier works.
	The normalization needed for the positive result is proved in
	Section~\ref{sec:global}; the intrinsic matrix extraction and Hermitian
	comparison needed for the obstruction are proved in
	Section~\ref{sec:transcendental}.
	
	\subsection{The main positive result: algebraic spectral descent}
	
	Let $K\subset\C$ be stable under actual complex conjugation $c$, put
	$k=K^c$, and assume $c|_K\ne1$, so that $[K:k]=2$.
	For a matrix $g$, let $E(g)$ be the set of all its characteristic roots
	in the fixed complex realization, and write $L_g=K(E(g))$.
	When $g$ has entries in $K$, this is its characteristic polynomial's
	splitting field. For a $c$-stable extension $L/K$, the norm to the fixed
	field is $N_{L/L^c}(x)=xc(x)$.
	
	Write $M^*=c(M)^t$. A Hermitian model is a nonsingular matrix $H=H^*$,
	with special unitary group
	\[
	\SU(H)(k)=\{g\in\SL_m(K):g^*Hg=H\}.
	\]
	The normalized forms considered here are
	\[
	H=\mathbb H\perp\langle1,q_3,\ldots,q_{m-1}\rangle,
	\qquad \mathbb H=\begin{pmatrix}0&1\\1&0\end{pmatrix},
	\qquad m\ge4,
	\]
	and the standard form is $J_m=\mathbb H\perp I_{m-2}$.
	Their signature $(m-1,1)$ refers to the fixed complex embedding;
	no positivity at every real embedding of $k$ is assumed.
	Standard unitary descent over $L$ means that one change of basis
	realizes the entire prescribed group in $\SU(J_m)(L^c)$.
	Thus the matrices have entries in $L$, while the unitary algebraic
	group is defined over $L^c$. It suffices to have a Hermitian similarity
	$P^*HP=aJ_m$, with $P\in\GL_m(L)$ and $a\in(L^c)^\times$,
	since scaling a form does not change its special unitary group.
	
	\paragraph{Simultaneous norms.}
	Suppose that $K$ is algebraic over $\Q$ and that $H$ is normalized as
	above. For every $k$-Zariski-dense subgroup $\Gamma\subset\SU(H)(k)$,
	there is one regular semisimple loxodromic element $\gamma\in\Gamma$
	and there are $x_j\in L^\times$, where $L=K(E(\gamma))$, such that
	\[
	x_jc(x_j)=q_j\qquad(3\le j\le m-1).
	\]
	These simultaneous identities give a single change of basis for the
	whole group by rescaling the diagonal coordinates of the residual form.
	The same element can be chosen so that, for every $n\ge1$,
	$\gamma^n$ remains regular semisimple and loxodromic and
	\[
	K(E(\gamma^n))=L.
	\]
	The witnesses $x_j$ are independent of $n$. Here regular semisimplicity
	means that the characteristic roots are distinct; loxodromy means that
	one root has modulus greater than one and its paired root has modulus
	less than one. Density always means Zariski density over the indicated
	field of definition.
	
	Theorem~\ref{thm:norm-numeric} proves the norm assertion over number fields,
	and Corollary~\ref{cor:norm-realizing-powers} supplies the choice of powers.
	Proposition~\ref{prop:algebraic-transfer} extends both conclusions to
	algebraic $K$ of possibly infinite degree, by placing the required data
	at a finite number-field stage. No finite-generation or arithmeticity
	hypothesis on $\Gamma$ is imposed.
	
	\paragraph{Standard unitary descent.}
	For an arbitrary initial Hermitian form over such an algebraic $K$,
	of signature $(m-1,1)$,
	Corollary~\ref{cor:two-spectrum-algebraic} first uses the spectrum of a
	loxodromic element $A$ to normalize over $K(E(A))$. Applying the norm
	theorem there gives
	\[
	L=K(E(A),E(\gamma))\subset F_{\rm eig}(\Gamma),
	\qquad [L:K]\le(m!)^2,
	\]
	where $F_{\rm eig}(\Gamma)=K(E(g):g\in\Gamma)$ is the total eigenvalue field
	over the original $K$. Thus two spectra suffice without assuming
	isotropy over $K$ or finite generation of $\Gamma$; the preliminary
	spectrum of $A$ is part of this construction.
	This is the principal descent conclusion. In the normalized case, the
	norm-realizing element is selected for the group and its coefficients;
	the earlier obstruction for a prescribed initial spectrum is compatible
	with this choice. Our result proves the simultaneous algebraic version
	of the norm formulation in~\cite[Conjecture~6.4]{CanoCastroGusevskii},
	which allows a different element for each coefficient.
	Proposition~\ref{prop:global-transfer} gives the imaginary quadratic
	case of signature $(3,1)$, with degree at most $576$, and proves
	Conjecture~1.5 of that paper without its finite-generation hypothesis.
	
	\subsection{The main obstruction: intrinsic transcendental fields}
	
	Theorem~\ref{thm:transcendental-obstruction}, together with
	Corollary~\ref{cor:no-total-unitary-descent}, gives the negative answer
	beyond the algebraic setting: adjoining all spectra need not produce
	a standard unitary realization. The theorem constructs a single discrete, convex cocompact,
	Zariski-dense free subgroup $\Gamma\subset\SU(3,1)$ of rank four.
	For every $A\in\Gamma\setminus\{1\}$, put
	\[
	K_A=\Q(\tr\Gamma,E(A)),\qquad k_A=K_A^c.
	\]
	Every normalized extraction over $K_A$ produces a coefficient $q_A$
	that remains a non-norm after adjoining any finite family of spectra,
	including the empty family, and over the total spectral field.
	The group is fixed first, then $A$ and its extraction, and finally the
	finite family of words. A valuation testing that family may depend on
	the family; the group and the extracted coefficient do not.
	
	Although the fields $K_A$ can depend on $A$, their total spectral
	extensions agree:
	\[
	L_{\rm spec}(\Gamma)=\Q(E(g):g\in\Gamma)
	=K_A(E(g):g\in\Gamma).
	\]
	Every trace is a sum of characteristic roots, and $E(A)$ is already
	among the adjoined spectra. Corollary~\ref{cor:no-total-unitary-descent} excludes a standard
	unitary realization of this group over the total field, even if the
	initial conjugating matrix is allowed to have arbitrary complex entries.
	
	This construction shows that an unrestricted transcendental extension
	of the spectral norm principle fails. It lies outside the number-field
	setting of Proposition~6.3 and the imaginary quadratic hypothesis of
	Conjecture~1.5 in~\cite{CanoCastroGusevskii}; it does not contradict
	the algebraic descent result.
	
	\subsection{A separate local splitting theorem}
	
	Theorem~\ref{thm:split-local-spectra} addresses the proposed local
	even-degree assertion. For every
	quadratic extension $F/k$ of nonarchimedean local fields of characteristic
	zero, it constructs a discrete,
	torsion-free, Zariski-dense free subgroup $\Lambda$ of a four-dimensional
	quasi-split unitary group for which
	\[
	E_F(g)=F\qquad(g\in\Lambda).
	\]
	Here $E_F(g)$ denotes the splitting field over $F$, including $F$ itself.
	Every individual splitting degree is therefore one, disproving
	Conjecture~1.4 of Cano--Castro--Gusevskii~\cite{CanoCastroGusevskii}.
	The form has local Witt index two, distinct from complex signature
	$(3,1)$. Thus the global conjecture holds although this local assertion
	fails: the positive proof uses simultaneous norm realization and does
	not impose a universal even-degree splitting condition on local groups.
	Lemma~\ref{lem:uniform-products} controls words of arbitrary length;
	Section~\ref{subsec:explicit-local-example} gives an explicit example
	over $\Q_3(i)$.
	
	\subsection{Fields of definition and proof ideas}
	
	Several field invariants enter this problem, with different roles.
	For a Zariski-dense subgroup of a semisimple group in characteristic zero,
	Vinberg~\cite[Theorem~1]{Vinberg1971} identifies the field generated by
	adjoint traces as the minimal field of definition of the adjoint
	representation. This field must be distinguished from the natural
	trace field and from a field over which the Hermitian form becomes
	standard. Hamilton--Reid~\cite[Theorem~1.3]{HamiltonReid2004} prove that
	the total eigenvalue field of a non-elementary subgroup of the
	time-orientation-preserving group $\mathrm O_o(n,1;\R)$, $n\ge2$, has
	infinite degree over $\Q$. Over finitely generated infinite fields,
	Pisolkar--Rajan~\cite{PisolkarRajan2016} relate the commensurability class
	of the compositum of generic splitting fields to the absolute
	Killing--Cartan type, subject to their stated Weyl-group exceptions.
	That absolute-type question is distinct from standardizing a fixed
	Hermitian representation over its spectral field.
	
	Generic tori provide another way to obtain spectral information.
	Over a finitely generated field $K$,
	Prasad--Rapinchuk~\cite[Theorem~1]{PrasadRapinchuk2017} prove that the
	$K$-generic elements of a Zariski-dense subgroup $\Gamma\subset G(R)$
	of a semisimple group form an open subset for congruences modulo ideals
	of finite index, and hence in the profinite topology, when $R\subset K$
	is a finitely generated ring. The subgroup itself need not be finitely
	generated. Their Theorem~2 gives a generic element of infinite order
	in every Zariski-dense subsemigroup containing an element of infinite
	order when the ambient group is connected and absolutely almost simple.
	Their local--global theorem for maximal subfield embeddings compatible
	with unitary involutions~\cite[Theorem~A(i)]{PrasadRapinchuk2010} addresses
	the embedding problem over global fields. Here, the spectrum must also
	remove the norm classes of specified Hermitian coefficients inside the
	prescribed group. This requires control of the Hermitian embedding
	data as well as the abstract splitting torus.
	
	The local construction lies in the Schottky setting over local fields
	developed by Benoist~\cite[\S\S2.2, 4.3, 5.1]{Benoist1997}.
	A uniform product estimate in the first three exterior powers separates
	all four eigenvalue moduli for every nontrivial word. It also supports
	the explicit fifteen-generator example.
	
	The number-field proof proceeds from the correlations among local
	representations to one word realizing all the required norms.
	A blockwise Goursat argument identifies these correlations through
	the adjoint trace field and the rational source
	\cite{MilneAG2017}. A semigroup construction combines proximality
	\cite{Benoist1997} with return in compact factors and access to local
	target neighborhoods \cite{Glockner2016}. The finite list of connecting
	elements is chosen before the return exponent. The resulting spectral
	behavior removes every local norm obstruction for the same quadratic
	extension, so the cyclic Hasse norm theorem applies
	\cite{MilneCFT2020}. The local tools are developed explicitly in the body.
	
	For the transcendental construction, a common interval of convex
	cocompact deformations supplies one group with the required word
	properties. The selection and geometric inputs come from
	\cite{Borel1983,Winkelmann2002,GuichardWienhard2012}.
	Simultaneous continuation of finite spectral composita gives an
	integer-valued valuation invariant under actual complex conjugation.
	Every norm has even value, while the extracted Hermitian coefficient
	has odd value. A hypothetical norm witness in the total field would
	already occur at a finite spectral level, where this parity is impossible.
	
	\subsection{Organization}
	
	Section~\ref{sec:local} proves the local counterexample.
	Section~\ref{sec:goursat} identifies the rational source and its correlated
	blocks, preparing the simultaneous norm theorem and the persistence
	under powers in Section~\ref{sec:numeric}.
	Section~\ref{sec:global} passes to algebraic fields and constructs the
	common change of basis, including normalization by a preliminary spectrum.
	Section~\ref{sec:transcendental} develops the intrinsic transcendental
	obstruction. All essential arguments and the data for the explicit local
	example are contained in the article. The accompanying supplement provides
	exact computational checks and additional examples illustrating the hypotheses
	\cite{SpectralFieldsComputationsZenodo}.

	\section{A local counterexample with completely split spectra}
	\label{sec:local}
	
	Let $k$ be a nonarchimedean local field of characteristic zero, and let
	$F/k$ be a quadratic extension. Write $\overline{x}$ for its nontrivial
	involution and $A^*=\overline A^{\,t}$. Fix the absolute value on $k$
	and its unique extension to $\overline F$. In this section, $E_F(A)$
	denotes the splitting field over $F$ of the characteristic polynomial
	of $A$; in particular, it contains $F$. We consider the group
	\[
	\mathsf J=\begin{pmatrix}
		0&0&0&1\\0&0&1&0\\0&1&0&0\\1&0&0&0
	\end{pmatrix},\qquad
	G(k)=\{A\in\operatorname{SL}_4(F):A^*\mathsf J A=\mathsf J\}.
	\]
	The form $\mathsf J$ has Witt index two. This local convention neither
	imposes nor represents the complex signature $(3,1)$.
	
	\begin{theorem}
		\label{thm:split-local-spectra}
		For every quadratic extension $F/k$ as above, there is a discrete,
		Zariski-dense subgroup $\Lambda\subset G(k)$ that is free of rank
		fifteen and such that every $\gamma\in\Lambda\setminus\{I\}$ has
		four eigenvalues in $F$ with pairwise distinct absolute values.
		In particular,
		\[
		\forall\gamma\in\Lambda,\qquad
		E_F(\gamma)=F,\qquad [E_F(\gamma):F]=1.
		\]
		The same conclusion holds for every four-dimensional quasi-split
		Hermitian form over $F/k$.
	\end{theorem}
	
	The proof begins with an estimate uniform in the length of a cyclically
	reduced word. Applied in the first three exterior powers, it separates
	all four eigenvalue moduli. Conjugation then gives the assertion for
	every nontrivial word, including inverses.
	
	\subsection{A uniform product estimate}
	
	On a finite-dimensional vector space over $F$, we use the sup norm
	relative to a fixed basis, together with the induced norms on matrices
	and linear functionals. Each exterior power carries the exterior basis.
	
	\begin{lemma}
		\label{lem:uniform-products}
		Let $\mathcal A$ be a finite alphabet with an involution
		$s\mapsto s^{-1}$. On a vector space $W$ of dimension at least two, let
		\[
		B_s=P_s+E_s,\qquad P_s=v_s\varphi_s,\qquad
		\|v_s\|=1,\quad \varphi_s(v_s)=1.
		\]
		Suppose that $M\ge1$, $0<\delta\le1$, and $\varepsilon\ge0$ satisfy
		\[
		\|\varphi_s\|\le M,\qquad \|E_s\|\le\varepsilon,
		\qquad |\varphi_s(v_t)|\ge\delta\quad(t\ne s^{-1}),
		\qquad q=\frac{M\varepsilon}{\delta^2}<1.
		\]
		For a cyclically reduced word $s_1\cdots s_n$ of length $n\ge1$,
		set $s_{n+1}=s_1$ and
		\[
		B=B_{s_1}\cdots B_{s_n},\qquad
		T=\prod_{j=1}^n\varphi_{s_j}(v_{s_{j+1}}).
		\]
		Then
		\begin{equation}\label{eq:uniform-trace}
			|\operatorname{tr}B|=|T|\ge\delta^n,
			\qquad \left\|\bigwedge^2B\right\|
			\le(M\varepsilon)^n<|T|^2.
		\end{equation}
		The matrix $B$ has a simple eigenvalue $\beta_1\in F$ whose absolute
		value is strictly larger than those of all the other eigenvalues.
		If $\beta_2$ has the second-largest absolute value, with eigenvalues
		counted with multiplicity, then
		\[
		|\beta_1|=|T|,\qquad
		\frac{|\beta_2|}{|\beta_1|}\le q^n.
		\]
	\end{lemma}
	
	\begin{proof}
		In the expansion of $\operatorname{tr}\prod_j(P_{s_j}+E_{s_j})$,
		the term containing no error factors is $T$, by the composition
		formula for rank-one operators. Consider a term containing $m\ge1$
		error factors and at least one projector. In the cyclic order of the
		word, the error factors form $b$ nonempty blocks, where $1\le b\le m$.
		A block of $\ell$ error factors between two projectors replaces
		$\ell+1$ transitions in $T$ by the coefficient
		\[
		\varphi_s(E_{t_1}\cdots E_{t_\ell}v_u),
		\qquad
		|\varphi_s(E_{t_1}\cdots E_{t_\ell}v_u)|
		\le M\varepsilon^\ell.
		\]
		All other transitions remain unchanged. Cyclic reduction ensures that
		every original transition is admissible, including the transition
		from the last letter to the first. Dividing the absolute value of
		the term by $|T|$ gives the bound
		\[
		M^b\varepsilon^m\delta^{-(m+b)}
		\le (M\varepsilon/\delta^2)^m=q^m<1.
		\]
		Here we used $b\le m$, $M\ge1$, and $\delta\le1$.
		The term in which every factor is an error satisfies
		\[
		\frac{|\operatorname{tr}(E_{s_1}\cdots E_{s_n})|}{|T|}
		\le(\varepsilon/\delta)^n\le q^n.
		\]
		In an ultrametric sum, the number of terms does not increase the
		maximum bound. Hence
		$|\operatorname{tr}B-T|\le q|T|<|T|$, proving the asserted
		equality of absolute values without any additional dependence on $n$.
		
		Since $P_s$ has rank one, all its $2\times2$ minors vanish. Each
		$2\times2$ minor of $P_s+E_s$ is bounded by
		$\max(M\varepsilon,\varepsilon^2)=M\varepsilon$, because
		$\varepsilon<\delta^2/M\le1\le M$. Functoriality of exterior powers
		and submultiplicativity of the norm therefore give
		\[
		\|\bigwedge^2 B\|\le(M\varepsilon)^n
		<\delta^{2n}\le |T|^2.
		\]
		
		Extend the coordinate sup norm to a finite splitting field. The absolute
		values of the matrix entries, and hence the induced matrix norm, are
		unchanged by this scalar extension. Let $r_1\ge r_2\ge\cdots$ be the
		absolute values of the eigenvalues, counted with algebraic
		multiplicity. The spectral-radius bound and the ultrametric
		inequality imply
		\[
		r_1r_2\le\|\bigwedge^2B\|,
		\qquad |\operatorname{tr}B|\le r_1.
		\]
		If $r_1=r_2$, these inequalities contradict
		\eqref{eq:uniform-trace}. Thus $r_1>r_2$: the dominant eigenvalue is
		unique and simple, and its absolute value equals that of the trace.
		It follows that
		$r_1=|T|$ and $r_2/r_1\le(M\varepsilon)^n/|T|^2\le q^n$.
		Finally, every automorphism of $\overline F/F$ preserves the extended
		absolute value and permutes the eigenvalues. It must therefore fix
		the unique dominant eigenvalue. Its minimal polynomial is separable
		in characteristic zero, so this eigenvalue belongs to $F$.
	\end{proof}
	
	\subsection{Conjugating elements and a common power}
	
	We next choose conjugates that are transverse in all three exterior
	powers and whose infinitesimal generators span the Lie algebra.
	A common power will then make the product estimate apply.
	
	Let $\pi$ be a uniformizer of $k$, and set
	\[
	\lambda(t)=\operatorname{diag}(t^{-3},t^{-1},t,t^3),\qquad
	a=\lambda(\pi),\qquad X=\operatorname{diag}(-3,-1,1,3).
	\]
	The homomorphism $\lambda:\mathbb G_m\to G$ is defined over $k$:
	the entries paired by $\mathsf J$ have product one, and the
	determinant is one. For $r=1,2,3$, write
	$\rho_r=\bigwedge_F^r$ and $(b_1,b_2,b_3)=(3,4,3)$.
	The dominant vectors of $\rho_r(a)$ and $\rho_r(a^{-1})$ are,
	respectively,
	\[
	v_{+,r}=e_1\wedge\cdots\wedge e_r,
	\qquad v_{-,r}=e_{5-r}\wedge\cdots\wedge e_4.
	\]
	Let $\varphi_{\pm,r}$ be the dual coordinate functionals, and put
	$P_{\pm,r}=v_{\pm,r}\varphi_{\pm,r}$. A diagonal calculation gives
	\begin{equation}\label{eq:diagonal-gap}
		\|\pi^{b_rN}\rho_r(a^{\pm N})-P_{\pm,r}\|
		\le |\pi|^{2N}\qquad(N\ge1).
	\end{equation}
	The gap between the first two exponents is two in all three
	exterior powers and for both signs.
	
	Choose fifteen matrices $h_i\in G(k)$ satisfying simultaneously
	\begin{align}
		&\{\operatorname{Ad}(h_i)X:1\le i\le15\}
		\text{ is a basis of }\mathfrak g=\operatorname{Lie}(G)(k),
		\label{eq:lie-basis}\\
		&\varphi_{\varepsilon,r}
		(\rho_r(h_i^{-1}h_j)v_{\eta,r})\ne0
		\quad(i\ne j,\ \varepsilon,\eta\in\{+,-\},\ 1\le r\le3).
		\label{eq:transversality}
	\end{align}
	We justify this choice before using it. The group $G$ splits over
	$F$ as $\operatorname{SL}_4$, so it is smooth and geometrically
	connected of dimension fifteen. Its Lie algebra $\mathfrak g$ is
	simple, since every ideal becomes an ideal of $\mathfrak{sl}_4(F)$
	after extending scalars. The subspace $V=\operatorname{span}_k\{\operatorname{Ad}(h)X:h\in G(k)\}$
	is invariant under $\operatorname{Ad}(G(k))$. Differentiating this
	invariance along $\operatorname{Ad}(\exp(tY))$, for sufficiently
	small $t$, gives $[Y,V]\subset V$. The limit remains in $V$
	because finite-dimensional subspaces are closed. Since
	$\mathfrak g$ is simple and $X\ne0$, we obtain $V=\mathfrak g$.
	There is therefore a family of fifteen conjugating elements
	satisfying \eqref{eq:lie-basis}, and this condition defines a
	nonempty Zariski-open subset of $G^{15}$.
	
	Each condition in \eqref{eq:transversality} also defines a
	nonempty open subset. For equal signs, nonemptiness follows by
	taking $h_i=h_j=I$; for opposite signs, take $h_i=I$ and
	$h_j=\mathsf J$. Indeed, $\mathsf J\in G(k)$ and it interchanges
	the extreme vectors up to sign. The coefficients are regular
	functions with values in $F$; taking their norms to $k$ expresses
	nonvanishing as an algebraic condition over $k$.
	
	The variety $G^{15}$ is geometrically integral. The preceding open
	subsets have nonempty intersection, and that intersection contains
	a $k$-rational point. To see the latter assertion, use the density
	of rational points in a smooth connected group over a local field
	of characteristic zero. An analytic chart around the identity
	contains an open ball of the same dimension as the group. A regular
	function vanishing at every rational point vanishes on that ball,
	so its convergent power series is zero. The injection of the
	regular local ring into its completion, together with the
	integrality of the group, then shows that the function is
	identically zero. The same argument applies to $G^{15}$.
	This proves the existence of the required family.
	
	For the alphabet $s=(i,\varepsilon)$, with
	$s^{-1}=(i,-\varepsilon)$, write
	\[
	P_{s,r}=\rho_r(h_i)P_{\varepsilon,r}\rho_r(h_i)^{-1}
	=v_{s,r}\varphi_{s,r}.
	\]
	Rescale each vector so that $\|v_{s,r}\|=1$, and rescale its
	functional inversely. Then $\varphi_{s,r}(v_{s,r})=1$.
	This normalization is possible because the sup norm of a nonzero
	vector belongs to the value group of $F^\times$.
	By \eqref{eq:transversality}, all transitions with
	$t\ne s^{-1}$ are nonzero; when $t=s$, their value is one.
	Taking maxima and minima over all indices, define
	\begin{align*}
		M&=\max\bigl(1,\|\varphi_{s,r}\|\bigr),\\
		\delta&=\min\bigl(1,|\varphi_{s,r}(v_{t,r})|:
		t\ne s^{-1}\bigr),\\
		C&=\max\bigl(1,\|\rho_r(h_i)\|\,\|\rho_r(h_i^{-1})\|\bigr).
	\end{align*}
	Thus $M,C\ge1$ and $0<\delta\le1$. After fixing the conjugating
	elements, choose a single exponent $N\ge1$ such that
	\begin{equation}\label{eq:common-local-power}
		q=MC|\pi|^{2N}\delta^{-2}<1,
		\qquad R=|\pi|^{-3N}\delta>1,
	\end{equation}
	and set
	\[
	g_i=h_i a^Nh_i^{-1},\qquad
	\Lambda=\langle g_1,\ldots,g_{15}\rangle.
	\]
	In each exterior power, the matrices
	$B_{(i,\pm),r}=\pi^{b_rN}\rho_r(g_i^{\pm1})$
	satisfy the lemma with error $\varepsilon=C|\pi|^{2N}$.
	The constants and the exponent $N$ are common to all letters, signs,
	and exterior powers; no further choice depends on the word.
	
	\subsection{Spectra, freeness, discreteness, and density}
	
	Let $w$ be a cyclically reduced word of length $n\ge1$.
	The normalized product in the $r$th exterior power is
	$\pi^{b_rNn}\rho_r(w)$. By the lemma, it has a unique eigenvalue
	of strictly largest absolute value. Let
	$\mu_1,\ldots,\mu_4$ be the eigenvalues of $w$, ordered by
	nonincreasing absolute value and counted with multiplicity.
	The dominant eigenvalue in the $r$th exterior power has absolute
	value $|\mu_1\cdots\mu_r|$. If
	$|\mu_r|=|\mu_{r+1}|$, that exterior power would have two
	dominant eigenvalues, counted with multiplicity.
	Applying the lemma for $r=1,2,3$ gives
	\begin{equation}\label{eq:local-separation}
		|\mu_1|>|\mu_2|>|\mu_3|>|\mu_4|,
		\qquad\frac{|\mu_{r+1}|}{|\mu_r|}\le q^n.
	\end{equation}
	Each eigenvalue has a distinct absolute value. Every automorphism
	of $\overline F/F$ preserves absolute values and hence fixes each
	eigenvalue. Thus all four eigenvalues belong to $F$.
	Moreover, the trace bound in the first exterior power gives
	\begin{equation}\label{eq:local-radius}
		|\mu_1|\ge |\pi|^{-3Nn}\delta^n=R^n>1.
	\end{equation}
	
	Every nonempty reduced formal word is formally conjugate to a
	nonempty cyclically reduced word. Matrix conjugation preserves
	the spectrum, so no nonempty reduced word evaluates to $I$.
	Consequently $\Lambda$ is free of rank fifteen and torsion-free.
	The same reduction proves the stated spectral properties for
	every nontrivial element of $\Lambda$. The characteristic
	polynomial of the identity also splits over $F$.
	
	For discreteness, consider the identity neighborhood
	$\mathcal U=\{I+A:\|A\|<1\}\subset\operatorname{GL}_4(F)$.
	All eigenvalues of $A$ have absolute value less than one, so all
	eigenvalues of $I+A$ have absolute value one. Every nontrivial
	element of $\Lambda$, however, has spectral radius at least
	$R>1$. Therefore $\Lambda\cap\mathcal U=\{I\}$.
	
	Finally, the powers of $\pi^N$ are Zariski-dense in
	$\mathbb G_m$, since a nonzero Laurent polynomial cannot vanish
	at infinitely many distinct points. If $H$ is the Zariski closure
	of $\Lambda$ in $G$, it contains the tori
	$h_i\lambda(\mathbb G_m)h_i^{-1}$. Hence its Lie algebra contains
	all the vectors $\operatorname{Ad}(h_i)X$ and equals
	$\mathfrak g$ by \eqref{eq:lie-basis}. Algebraic groups in
	characteristic zero are smooth. Equality of dimensions and the
	connectedness of $G$ therefore imply $H=G$.
	
	A four-dimensional quasi-split Hermitian form has Witt index two
	and is the orthogonal sum of two hyperbolic planes. In a suitable
	basis, its matrix is therefore $\mathsf J$. Changing basis
	transports the constructed subgroup and preserves all the stated
	properties. This completes the proof of
	Theorem~\ref{thm:split-local-spectra}.
	
	\subsection{An explicit certificate over \texorpdfstring{$\mathbb Q_3(i)$}{Q3(i)}}
	\label{subsec:explicit-local-example}
	
	The preceding existence proof does not depend on the calculation
	below, which provides a counterexample with fully specified
	generators. Take $k=\mathbb Q_3$, $F=\mathbb Q_3(i)$,
	$i^2=-1$, and $|3|=3^{-1}$. This is an unramified quadratic
	extension, since $x^2+1$ is irreducible modulo three.
	For integers $a,b,c,d,u$, set $z_1=a+bi$, $z_2=c+di$, and
	\[
	N(a,b,c,d,u)=\begin{pmatrix}
		1&0&0&0\\z_1&1&0&0\\z_2&0&1&0\\
		-(ac+bd)+ui&-\overline z_2&-\overline z_1&1
	\end{pmatrix}.
	\]
	With $w=-(ac+bd)+ui$, the isotropy condition for the first column is
	\[
	w+\overline w+\overline z_1z_2+\overline z_2z_1=0.
	\]
	The other column pairings give $N^*\mathsf JN=\mathsf J$; triangularity
	gives $\det N=1$. These identities hold for every choice of the five
	integer parameters. The following five-tuples specify
	$h_j=N(p_j)\mathsf JN(q_j)\mathsf J$.
	\begin{center}
		\small
		\begin{tabular}{r|rrrrr|rrrrr}
			$j$&\multicolumn{5}{c|}{$p_j$}&\multicolumn{5}{c}{$q_j$}\\
			\hline
			1&0&-2&1&-1&-2&1&2&2&2&2\\
			2&2&0&2&-1&1&2&1&2&2&-2\\
			3&-1&-1&-1&0&0&1&2&-1&-2&0\\
			4&-1&-2&2&-2&0&0&1&2&0&1\\
			5&0&-2&-1&2&-1&2&0&1&1&-2\\
			6&0&0&0&-1&0&1&1&-2&-1&1\\
			7&-2&0&1&1&0&-1&-2&-2&1&-1\\
			8&-1&2&0&2&2&-2&-2&-2&-1&-2\\
			9&-2&2&-2&2&1&1&0&-1&2&0\\
			10&1&0&1&-2&1&1&0&1&2&0\\
			11&-2&-2&2&1&-1&-1&-1&2&2&-1\\
			12&-1&2&-1&1&2&-2&0&2&0&1\\
			13&1&2&1&0&1&1&2&1&1&2\\
			14&1&2&1&2&1&1&0&-2&-1&-1\\
			15&0&1&-2&0&1&-1&1&2&0&1\\
			\hline
		\end{tabular}
	\end{center}
	The generators are
	\[
	g_j=h_j\operatorname{diag}(3^{-18},3^{-6},3^6,3^{18})h_j^{-1}.
	\]
	Let $I_{+,r}=\{1,\ldots,r\}$ and
	$I_{-,r}=\{5-r,\ldots,4\}$. The transition coefficients are
	\[
	c_{j,l,\varepsilon,\eta,r}
	=\det((h_j^{-1}h_l)_{I_{\varepsilon,r},I_{\eta,r}}),
	\qquad j\ne l,\quad r=1,2,3.
	\]
	All $15\cdot14\cdot4\cdot3=2520$ coefficients are nonzero.
	With the normalization $v_3(3)=1$, the valuations
	$0,1,2,3,4,5$ occur with respective multiplicities
	$2054,360,94,6,4,2$.
	For the Lie-algebra calculation, set $K_j=\mathsf Jh_jXh_j^{-1}$.
	The matrix $K_j$ is skew-Hermitian. We use the following fifteen
	coordinates, in the displayed order:
	\[
	\begin{split}
		&\Im K_{11},\Im K_{22},\Im K_{33},\Im K_{44};\\
		&\Re K_{12},\Re K_{13},\Re K_{14},
		\Re K_{23},\Re K_{24},\Re K_{34};\\
		&\Im K_{12},\Im K_{13},\Im K_{14},\Im K_{24},\Im K_{34}.
	\end{split}
	\]
	Here $\Re$ and $\Im$ denote the coordinates in the basis $(1,i)$
	of $F/k$. The trace condition gives
	$\Im K_{23}=-\Im K_{14}$, explaining the omitted coordinate.
	The determinant of the fifteen rows is
	\[
	2278188751943281589354496\equiv1\pmod5,
	\]
	which proves linear independence over $\mathbb Q_3$.
	
	The matrices $h_j$ and $h_j^{-1}$ are integral, as are their
	exterior powers. They act isometrically for the sup norm, so
	the conjugated extreme vectors and their coordinate functionals
	have norm one. We may therefore take
	\[
	C=M=1,\quad\delta=3^{-5},\quad N=6,\quad
	\varepsilon=3^{-12},\quad q=\frac19,\quad R=3^{13}.
	\]
	Lemma~\ref{lem:uniform-products} then gives, for every cyclically
	reduced word of length $n\ge1$, a spectral radius of at least
	$3^{13n}$ and successive ratios of absolute values at most
	$9^{-n}$. The finite data above can be reproduced using Gaussian
	integer arithmetic and, as a second calculation,
	the corresponding $8\times8$ real matrix representation with
	integral entries and Bareiss elimination. The computational
	supplement contains the program and its output
	\cite{SpectralFieldsComputationsZenodo}. The assertion
	for words of arbitrary length follows from the proof above.
	
	\begin{remark}
		\label{rem:local-and-global}
		For a unitary matrix, the multiset of absolute values of its
		eigenvalues is invariant under $r\mapsto r^{-1}$. The four
		distinct absolute values in this counterexample are therefore
		$r_1,r_2,r_2^{-1},r_1^{-1}$, where $r_1>r_2>1$.
		This uses local Witt index two. It does not give a counterexample
		in complex signature $(3,1)$: the explicit diagonal matrix has
		two expanding eigenlines and cannot preserve a complex Hermitian
		form of Witt index one. Indeed, for real eigenvalues
		$\alpha,\beta>1$ and corresponding eigenvectors $v,w$, the identity
		$h(v,w)=\alpha\beta h(v,w)$ forces their span to be totally
		isotropic. Thus the local conclusion is separate from the
		norm and global descent questions.
	\end{remark}

	\section{Goursat blocks and the rational source}
	\label{sec:goursat}
	
	In the simultaneous norm argument, local representations arising from
	different embeddings need not vary independently. We identify these
	correlations through the algebraic group determined by adjoint traces.
	The block isomorphisms obtained here are algebraic; openness of local
	subgroups will be proved separately.
	
	\begin{lemma}
		\label{lem:goursat-blocks}
		Let $\Omega$ be an algebraically closed field of characteristic zero,
		and let $Y\subset\prod_{i\in I}A_i$ be a connected algebraic subgroup
		that projects surjectively onto each simple adjoint group $A_i$.
		There is a partition $I=\coprod_{B\in\mathcal B}B$ such that
		\[
		Y=\prod_{B\in\mathcal B}D_B,
		\qquad D_B=\{(\phi_i(a))_{i\in B}:a\in A_{i_B}\},
		\]
		where $i_B\in B$, $\phi_{i_B}=\mathrm{id}$, and
		$\phi_i:A_{i_B}\longrightarrow A_i$ are algebraic isomorphisms.
		For $i\in B$, the restriction
		\[
		p_i|_{D_B}:D_B\longrightarrow A_i
		\]
		is an isomorphism, whereas
		\begin{equation}
			\ker(p_i:Y\longrightarrow A_i)
			=\prod_{B'\ne B}D_{B'}.
			\label{eq:goursat-kernel}
		\end{equation}
		In particular, $p_i$ is an isomorphism on all of $Y$ if and only
		if $\mathcal B$ consists of a single block.
		
		If $G_i\longrightarrow A_i$ are the simply connected covers and
		$Z\subset\prod_iG_i$ is connected and subdirect, the same conclusion
		holds for $Z$, with blocks $\widetilde D_B$ defined by the unique
		lifts of the $\phi_i$. In this case,
		$p_i|_{\widetilde D_B}:\widetilde D_B\to G_i$ is an isomorphism
		with trivial kernel.
	\end{lemma}
	
	\begin{proof}
		For two factors, set $N_1=\{a:(a,1)\in Y\}$ and
		$N_2=\{b:(1,b)\in Y\}$. Surjectivity of the projections identifies
		$A_1/N_1$ with $A_2/N_2$. The subgroups $N_i$ are normal.
		A simple adjoint group has no proper nontrivial normal algebraic
		subgroup: a finite normal subgroup would be central, and a
		positive-dimensional one would contradict simplicity. Thus either
		$Y=A_1\times A_2$ or $Y$ is the graph of an isomorphism.
		
		For the induction step, add one factor to a product of blocks.
		If the new factor is not independent, it is a simple quotient of
		that product. The images of distinct blocks commute, so at most
		one can be the entire new, noncommutative simple group. The new
		factor is therefore linked to exactly one block by an isomorphism.
		This gives the asserted partition. The explicit description of
		each graph proves \eqref{eq:goursat-kernel}.
		
		For the simply connected version, let
		$q:\prod_iG_i\to\prod_iA_i$. Each $\phi_i$ lifts uniquely to the
		simply connected covers; see \cite[24.9]{MilneAG2017}. The lifted
		graph has the same dimension as $D_B$, and $q^{-1}(D_B)$ is a
		finite union of its central translates. Its identity component is
		therefore that graph. The image of $Z$ is $Y$, the kernel of $q$
		is finite, and $Z$ is connected. Comparing dimensions gives
		$Z=(q^{-1}(Y))^\circ=\prod_B\widetilde D_B$. The kernel formula
		and the isomorphism on each block follow from this product
		decomposition.
	\end{proof}
	
	\begin{proposition}
		\label{prop:rational-source}
		Let $k/\mathbb Q$ be a number field, let $G/k$ be absolutely
		almost simple and simply connected, and let $\Gamma\subset G(k)$
		be a $k$-Zariski-dense subgroup. Define
		\[
		R=\operatorname{Res}_{k/\mathbb Q}G,\quad
		Z=\overline\Gamma^{\mathrm{Zar},\mathbb Q}\subset R,\quad
		f=\mathbb Q(\operatorname{tr}_{k}\operatorname{Ad}\gamma:
		\gamma\in\Gamma).
		\]
		Here $\operatorname{tr}_{k}\operatorname{Ad}$ denotes the trace of
		the adjoint representation as a $k$-linear map; the field trace
		$\operatorname{Tr}_{k/\mathbb Q}$ is not applied.
		There is an absolutely almost simple, simply connected group
		$G_0/f$ such that
		\begin{equation}
			Z^\circ\simeq\operatorname{Res}_{f/\mathbb Q}G_0,
			\qquad (G_0)_k\simeq G.
			\label{eq:rational-source}
		\end{equation}
		The subgroup $\Gamma_0=\Gamma\cap Z^\circ(\mathbb Q)$ has finite
		index and is $\mathbb Q$-Zariski-dense in $Z^\circ$.
		The geometric blocks are precisely the fibers of
		\[
		\operatorname{Hom}(k,\overline{\mathbb Q})
		\longrightarrow\operatorname{Hom}(f,\overline{\mathbb Q}).
		\]
	\end{proposition}
	
	\begin{proof}
		The traces belong to $k$, so $f$ is a finite intermediate
		extension. Since $Z/Z^\circ$ is finite and $\Gamma$ is dense
		in $Z$, the subgroup $\Gamma_0$ has finite index and is dense
		in $Z^\circ$.
		Over $\overline{\mathbb Q}$, every projection of $Z^\circ$ onto
		a factor $G^\sigma$ is surjective. Indeed, the projection of $Z$
		is surjective by density, and a finite union of translates of a
		closed subgroup of smaller dimension cannot cover the connected
		group $G^\sigma$.
		
		The adjoint image $Y$ of $Z^\circ$ decomposes as in
		Lemma~\ref{lem:goursat-blocks}. A graph of an isomorphism
		between adjoint groups is its own normalizer: if $(a,b)$
		normalizes the graph of $\phi$, then $\phi(a)^{-1}b$ centralizes
		the entire second factor and hence is the identity. The same
		holds for a product of graphs, since conjugation in the ambient
		product does not permute its coordinates. As $Z$ normalizes
		$Z^\circ$, the adjoint image of every element of $\Gamma$
		belongs to $Y$.
		
		Within a block, the differential of the graph isomorphism
		intertwines the adjoint representations, including when the
		isomorphism is outer. Thus two embeddings in the same block
		agree on all adjoint traces and hence on $f$. For distinct
		blocks, the difference of the two adjoint characters is not
		identically zero on $Y$: one block can vary while the other
		remains at the identity. The adjoint character is nonconstant,
		as its root decomposition on a torus shows. Density of
		$\Gamma_0$ therefore supplies an element whose traces
		distinguish the two embeddings. Consequently, two embeddings
		belong to the same block if and only if they agree on $f$.
		
		The lifted blocks are simply connected, and Galois permutes
		them as it permutes the embeddings of $f$. Descent of an orbit
		of simply connected factors gives a restriction of scalars
		by \cite[Theorem 24.3]{MilneAG2017}. Here the orbit is
		transitive and its stabilizer has fixed field $f$. This proves
		the first isomorphism in \eqref{eq:rational-source}.
		
		Let $\sigma_0:k\hookrightarrow\overline{\mathbb Q}$ be the
		original embedding and let $B_0$ be its block. The restricted projection $\widetilde D_{B_0}\longrightarrow G^{\sigma_0}$
		is an isomorphism equivariant under
		$\operatorname{Gal}(\overline{\mathbb Q}/k)$. It descends to
		$(G_0)_k\simeq G$. The projection from all of
		$Z^\circ_{\overline{\mathbb Q}}$, however, has the kernel in
		\eqref{eq:goursat-kernel}. No step requires finite generation.
	\end{proof}
	
	\paragraph{Local factors and transported representations.}
	The dimension of $Z^\circ$ is $[f:\mathbb Q]\dim G$.
	When $[f:\mathbb Q]>1$, the projection from the whole group
	to one geometric factor cannot be an isomorphism.
	After descent,
	\begin{equation}
		(Z^\circ)_{\mathbb Q_p}
		\simeq\prod_{w\mid p}\operatorname{Res}_{f_w/\mathbb Q_p}(G_0)_{f_w}.
		\label{eq:source-product}
	\end{equation}
	The approximation variables are indexed by the places $w$ of $f$,
	each occurring once. The places $v\mid w$ of $k$ give correlated
	realizations of the same source factor. In type $A$, an
	outer isomorphism may replace the natural representation by its
	dual. Adjoint traces therefore identify the blocks, while the
	natural representation must be transported through the specified
	isomorphism. Duality inverts the eigenvalues and preserves their
	splitting field. It does not imply equality of natural traces
	or characteristic polynomials in two realizations related by
	an outer isomorphism.
	
	\section{Simultaneous norm realization over number fields}
	\label{sec:numeric}
	
	\begin{theorem}
		\label{thm:norm-numeric}
		Let $K\subset\mathbb C$ be a number field stable under actual
		complex conjugation $c$. Put $k=K^c$ and assume $[K:k]=2$.
		Let $m\ge4$, and let
		\[
		H=\mathbb H\perp\langle1,q_3,\ldots,q_{m-1}\rangle,
		\quad q_j\in k^\times,
		\qquad \mathbb H=\begin{pmatrix}0&1\\1&0\end{pmatrix},
		\]
		have signature $(m-1,1)$ at the distinguished embedding.
		For every $k$-Zariski-dense subgroup
		$\Gamma\subset\mathrm{SU}(H)(k)$, there is a regular semisimple
		loxodromic element $\gamma\in\Gamma$ such that
		\begin{equation}
			L=K(E(\gamma)),\qquad
			q_j\in N_{L/L^c}(L^\times)\quad(3\le j\le m-1).
			\label{eq:simultaneous-norms}
		\end{equation}
		Here $E(\gamma)$ denotes the set of all eigenvalues of the original matrix.
		One word is chosen after $\Gamma$ and all the $q_j$ have been
		fixed. No assumption of finite generation, arithmeticity, local
		openness of $\Gamma$, total positivity of the coefficients, or
		equality of $k$ with the adjoint trace field is required.
	\end{theorem}
	
	\noindent\textit{Outline of the proof.}
	Only finitely many places obstruct the norm identities. The rational
	source identifies correlated local representations, and
	Proposition~\ref{prop:filtered-expansion} imposes expansion on every
	word of one dense semigroup. For the residual factors,
	Propositions~\ref{prop:return} and~\ref{prop:effective-compact} give
	access to prescribed spectra. Corollary~\ref{cor:one-word} chooses
	connectors before a common return exponent;
	Proposition~\ref{prop:norm-assembly} then gives the simultaneous norms.
	The proof follows the construction of the local targets;
	Section~\ref{sec:global} extends the result to algebraic fields.
	
	\subsection{Local source groups and their representations}
	\label{subsec:local-source}
	We first identify which local representations come from the same
	source, so that their spectral conditions can be imposed together.
	Fix $\Gamma_0$, $f$, $G_0$, and
	$\mathcal Z=\operatorname{Res}_{f/\mathbb Q}G_0$ as in
	Proposition~\ref{prop:rational-source}. The Lie algebra of $\mathcal Z$
	is simple over $\mathbb Q$: over an algebraic closure, Galois
	permutes its simple ideals transitively, whereas a rational
	ideal would give a Galois-stable sum of such factors.
	
	\begin{lemma}
		\label{lem:source-models}
		Let $v$ be a finite place of $k$ for which $K_v/k_v$ is a
		field extension, and put $w=v|_f$. There are a quadratic
		extension $E_w/F_w$, with $F_w=f_w$, and an $m$-dimensional
		Hermitian form $h_w$ such that
		\[
		(G_0)_{F_w}\simeq\operatorname{SU}(h_w),\qquad
		E_w\otimes_{F_w}k_v\simeq K_v.
		\]
		After extension to $K_v$, the original representation is
		conjugate to the natural representation of this source group
		or to its dual. The natural representation and its exterior
		powers of orders $1,\ldots,m-1$ are absolutely irreducible
		over $E_w$.
	\end{lemma}
	\begin{proof}
		An inner form remains inner after a field extension, since a
		trivial Galois action on the Dynkin diagram remains trivial
		upon restriction. The source group is therefore outer.
		The classification of type $A$
		\cite[Theorem 24.44 and Remark 24.45]{MilneAG2017}
		describes it by a central simple algebra $D/E_w$ with a
		unitary involution $*$. Its center encodes the diagram
		action and becomes the quadratic extension $K_v/k_v$.
		
		We prove the required corestriction statement. For a
		quadratic extension $E/F$ and an algebra $D/E$ with a
		unitary involution, let $ {}^\sigma D$ denote the conjugate
		algebra. The isomorphism
		\[
		D\otimes_E{}^\sigma D\longrightarrow\operatorname{End}_E(D),
		\qquad a\otimes{}^\sigma b\longmapsto(x\mapsto axb^*)
		\]
		identifies the semilinear interchange of the two factors
		with $T\mapsto JTJ$, where $J(x)=x^*$.
		The fixed space $D^J$ descends the underlying $E$-vector
		space $D$ to $F$: choosing nonzero $u$ with $\sigma(u)=-u$,
		we have
		\[
		x=\frac{x+x^*}{2}+u\frac{x-x^*}{2u},
		\]
		and both coefficients belong to $D^J$. The fixed algebra
		defining the corestriction is therefore
		$\operatorname{End}_F(D^J)$, which is split. For local fields,
		\cite[III, equation (30)]{MilneCFT2020} gives
		\[
		\operatorname{inv}_F\operatorname{Cor}_{E/F}[D]
		=\operatorname{inv}_E[D].
		\]
		Injectivity of the invariant map forces $[D]=0$, so
		$D\simeq M_m(E_w)$. A unitary involution on a matrix algebra
		is adjoint to a Hermitian form: writing it as
		$X\mapsto Q^{-1}\overline X^{\,t}Q$, the condition that its
		square be the identity gives $\overline Q^{\,t}=aQ$ with
		$a\overline a=1$, and Hilbert 90 allows us to rescale $Q$.
		
		Every automorphism of the split group $\operatorname{SL}_m$
		is inner or is the composite of an inner automorphism with
		inverse transpose \cite[\S24f]{MilneAG2017}. This proves
		the assertion about the natural representation and its dual.
		The specified simply connected isomorphism introduces no
		variable scalar character. The source group splits over
		$E_w$, and its exterior powers are the irreducible exterior
		powers of the standard representation of
		$\operatorname{SL}_m$; the other factors of the restriction
		of scalars act trivially.
		After transport, both Hermitian forms are invariant under
		the same irreducible representation. If their matrices are
		$H_1,H_2$, then $H_1^{-1}H_2$ commutes with the representation
		and is scalar by Schur's lemma. Hermitian symmetry places
		this scalar in $k_v^\times$. Thus the original form is
		similar to the transported source form, or to its conjugate
		or dual, according to the chosen isomorphism.
	\end{proof}
	
	Let $P$ be the union of the finite places at which some $q_j$
	is not a norm, and let $W$ be the set of their restrictions
	to $f$. These sets are finite: away from ramification, dyadic
	places, and zeros and poles of the coefficients, units are
	norms in an unramified quadratic extension. Each source
	place $w$ occurs once, even when several $v$ lie above it;
	these realizations are not independent variables.
	
	Suppose that $v$ is real, $K_v=\mathbb C$, and some
	$q_j(v)<0$. The summands $\mathbb H$, $\langle1\rangle$,
	and $\langle q_j\rangle$ provide at least two positive and
	two negative directions. Here $f_w=k_v=\mathbb R$.
	Models above the same $w$ have the same signature up to
	interchanging its two entries, and duality preserves the
	numbers of expanding and contracting eigenvalues.
	These conditions are compatible: no such source of real
	rank at least two is the distinguished source of real
	rank one.
	
	\subsection{Semigroup filtering and simultaneous expansion}
	
	We shall impose the expanding spectral conditions on every word of a
	single dense semigroup. Passing to a subsemigroup then preserves each
	condition already obtained. The remaining finite places will be treated
	through the compact blocks of one common expanding element.
	
	\begin{lemma}
		\label{lem:proximal-semigroup}
		A Zariski-dense semigroup in a connected real semisimple
		group contains an element whose Jordan projection lies
		in the interior of a restricted Weyl chamber.
		In $\operatorname{SU}(p,q)$, $p\ge q$, the moduli in its
		natural complex representation are
		$e^{u_1},\ldots,e^{u_q},1,\ldots,1,e^{-u_q},\ldots,e^{-u_1}$,
		where $u_1>\cdots>u_q>0$.
		Its exterior powers of orders $j\le q$ have a simple
		dominant eigenvalue, and their normalized powers converge
		to the corresponding rank-one projector.
	\end{lemma}
	\begin{proof}
		The simultaneous proximality result in
		\cite[\S3.1]{Benoist1997}, together with
		$\theta_\Gamma=\Pi$ over the real numbers and the
		proposition in \cite[\S3.2]{Benoist1997}, proves the first
		assertion for semigroups. The weights of a maximal
		real split torus give the displayed moduli.
		In the $j$th exterior power, with $j\le q$, the product
		of the $j$ largest expanding eigenvalues is the unique
		eigenvalue of maximal modulus and has multiplicity one.
		All other Jordan blocks decay after normalization of
		the powers. The conclusion uses these proximality
		results, not merely the fact that the limit cone has
		nonempty interior, and does not require the initial
		element to be algebraically regular.
	\end{proof}
	
	\begin{lemma}
		\label{lem:wedge-escape}
		Let $g$ be a unitary matrix of dimension $m$ over a
		quadratic extension of local fields. If
		$\bigwedge^r g$ is proximal and $1\le r\le m/2$, then
		$g$ has at least $r$ eigenvalues of modulus greater
		than one.
	\end{lemma}
	\begin{proof}
		List the moduli with multiplicity as
		$r_1\ge\cdots\ge r_m$. Reciprocity gives
		$r_{m+1-j}=r_j^{-1}$, and proximality of the exterior
		power gives $r_r>r_{r+1}$. Symmetry implies $r_r\ge1$.
		If $r_r=1$, then $r_{m+1-r}=1$, and monotonicity would
		force $r_{r+1}=1$. Hence $r_r>1$. The argument applies
		to both real unitary factors and nonarchimedean ones.
	\end{proof}
	
	\begin{lemma}
		\label{lem:compact-return}
		Let $E/\mathbb Q_p$ be finite, let $A\in\GL_d(E)$, and
		let $E'/E$ be a finite extension that splits its characteristic
		polynomial. Choose an eigenvalue $\beta\in E'$ of maximal
		modulus. Let $P_{\max}\in\End_{E'}(E'^d)$ be the spectral
		projector onto the sum of the generalized eigenspaces
		corresponding to eigenvalues $\lambda$ with
		$|\lambda|=|\beta|$. Then
		\begin{equation}\label{eq:compact-return}
			\lim_{n\to\infty}\beta^{-n!}A^{n!}=P_{\max}
			\qquad\text{in }\End_{E'}(E'^d).
		\end{equation}
		The sequence $n!$ works simultaneously for any finite
		collection of matrices and local fields, including extensions
		of the same source completion.
	\end{lemma}
	\begin{proof}
		The group $\mathcal O_{E'}^\times$ is profinite. For a
		unit $u$ and any finite quotient of this group, the order
		of the image of $u$ divides $n!$ for all sufficiently
		large $n$. Hence $u^{n!}\to1$. Moreover, for each fixed
		integer $j\ge1$,
		\[
		\binom{n!}{j}=\frac{n!(n!-1)\cdots(n!-j+1)}{j!}
		\longrightarrow0\quad\text{in }\mathbb Q_p,
		\]
		by continuity of the polynomial $x\mapsto\binom{x}{j}$
		and the convergence $n!\to0$ in $\mathbb Q_p$.
		
		In a Jordan decomposition over $E'$, every block with
		eigenvalue $\lambda$ has the form $\lambda(I+N)$, where
		$N$ is nilpotent after rescaling the nilpotent part.
		Thus
		\[
		\beta^{-n!}\bigl(\lambda(I+N)\bigr)^{n!}
		= (\lambda/\beta)^{n!}
		\sum_{j=0}^{d-1}\binom{n!}{j}N^j.
		\]
		If $|\lambda|=|\beta|$, the first factor tends to one
		and the terms with $j\ge1$ tend to zero. If
		$|\lambda|<|\beta|$, the first factor tends to zero
		and the sum remains bounded. The limit is the
		identity on each dominant generalized eigenspace
		and zero on the others. It is therefore $P_{\max}$.
		Neither semisimplicity of $A$ nor uniqueness of the
		dominant eigenvalue is needed.
	\end{proof}
	
	\begin{lemma}
		\label{lem:tracebound}
		Let $E/\mathbb Q_p$ be finite, and let $T\subset\SL_d(E)$
		be a semigroup whose $E$-linear span is $M_d(E)$.
		If $\{|\operatorname{tr} t|_E:t\in T\}$ is bounded,
		then $T$ is relatively compact in $\SL_d(E)$.
		In particular, if $T$ is unbounded, some $t\in T$
		satisfies $|\operatorname{tr} t|_E>1$.
	\end{lemma}
	\begin{proof}
		Choose $s_1,\ldots,s_{d^2}\in T$ forming a basis of
		$M_d(E)$. For every $t\in T$, the products $ts_j\in T$
		have uniformly bounded traces. The pairing
		$(X,Y)\mapsto\operatorname{tr}(XY)$ is nondegenerate
		on $M_d(E)$, so the $d^2$ linear forms
		$X\mapsto\operatorname{tr}(Xs_j)$ form a basis of
		the dual space. A fixed change of coordinates
		therefore bounds every entry of $t$ uniformly.
		Since $\det t=1$, the entries of $t^{-1}$ are
		polynomials in those of $t$ and are also uniformly
		bounded. The closure remains in $\SL_d(E)$ and
		is compact. For the last assertion, take the
		contrapositive: if no $t$ has
		$|\operatorname{tr}t|_E>1$, all traces are bounded
		by one.
	\end{proof}
	
	\begin{lemma}
		\label{lem:local-coefficients}
		Let $X/\mathbb Q$ be geometrically irreducible, and let
		$S\subset X(\mathbb Q)$ be Zariski-dense. Suppose that
		$f_1,\ldots,f_s$ are nonzero regular functions after scalar extension
		to a finite list of fields containing $\mathbb Q$. Then
		$S\cap\bigcap_i\{f_i\ne0\}$ is Zariski-dense, including after scalar
		extension. Rational functions are also allowed, provided that their
		denominators are included among the nonvanishing conditions.
	\end{lemma}
	\begin{proof}
		Write $f_i=\sum_j a_{ij}f_{ij}$, where
		$f_{ij}\in\mathbb Q[X]$ and the coefficients $a_{ij}$ are linearly
		independent over $\mathbb Q$. This expression involves only the
		linear span of the finitely many coefficients of $f_i$. At a rational
		point, $f_i$ vanishes if and only if every $f_{ij}$ vanishes. Choose
		one nonzero $f_{ij}$ for each $i$. Their product is nonzero in the
		domain $\mathbb Q[X]$. Its principal open set meets $S$ in a dense
		subset and is contained in the required intersection. The same
		separation of coefficients proves that this rational subset remains
		dense after scalar extension.
	\end{proof}
	
	\begin{lemma}
		\label{lem:common-trace}
		Let $\mathcal Z/\mathbb Q$ be a connected, geometrically
		irreducible algebraic group, and let
		$S\subset\mathcal Z(\mathbb Q)$ be a Zariski-dense
		semigroup. Fix a finite index set $I$, local fields
		$E_i/\mathbb Q_{p_i}$, and absolutely irreducible
		algebraic representations
		\[
		\rho_i:\mathcal Z_{E_i}\longrightarrow\SL_{d_i,E_i}
		\qquad(i\in I).
		\]
		Suppose that $\rho_i(S)$ is not relatively compact
		for any $i\in I$. Then, for every nonempty Zariski
		open subset $\mathcal U\subset\mathcal Z$ defined
		over $\mathbb Q$, there is
		$s\in S\cap\mathcal U(\mathbb Q)$ such that
		\begin{equation}\label{eq:all-traces}
			|\operatorname{tr}\rho_i(s)|_{E_i}>1
			\qquad\text{for all }i\in I.
		\end{equation}
		More precisely, the traces are simultaneously
		unbounded on $S\cap\mathcal U$: there is a
		sequence $s_l$ in this set such that
		$|\operatorname{tr}\rho_i(s_l)|_{E_i}\to\infty$
		for every $i$. No topological independence of
		the places, openness of $S$, or strong
		approximation is assumed.
	\end{lemma}
	\begin{proof}
		If $I$ is empty, Zariski density of $S$ suffices.
		Assume $I\ne\varnothing$. Absolute irreducibility
		and Burnside's theorem give
		$\operatorname{span}_{E_i}\rho_i(S)=M_{d_i}(E_i)$.
		Equivalently, this follows from the double
		centralizer theorem after extension to an
		algebraic closure and descent of the dimension;
		see \cite[IV, Theorem 1.14]{MilneCFT2020}.
		By Lemma~\ref{lem:tracebound}, for each $i$
		choose $g_i\in S$ with
		$|\operatorname{tr}\rho_i(g_i)|_{E_i}>1$.
		List these elements as $g_1,\ldots,g_b$, allowing
		repetitions. In every local representation
		$\rho_i$, the spectral radius of every
		$\rho_i(g_j)$ is at least one because its
		determinant is one. The spectral radius of
		$\rho_i(g_i)$ is greater than one by the
		ultrametric inequality for the trace.
		For each pair $(j,i)$, choose an eigenvalue
		$\beta_{j,i}$ of maximal modulus. Then
		\begin{equation}\label{eq:beta-growth}
			M_i:=\prod_{j=1}^b|\beta_{j,i}|_{E'_i}>1.
		\end{equation}
		Here $E'_i/E_i$ is a finite local extension
		splitting the finitely many spectra, with the
		absolute value extended from $E_i$.
		
		Let $P_{j,i}$ be the dominant projector supplied
		by Lemma~\ref{lem:compact-return} for $\rho_i(g_j)$.
		The multilinear expression
		\[
		F_i(h_1,\ldots,h_b)=
		\operatorname{tr}\bigl(P_{1,i}\rho_i(h_1)P_{2,i}\rho_i(h_2)
		\cdots P_{b,i}\rho_i(h_b)\bigr)
		\]
		is not identically zero on $\mathcal Z^b$.
		Indeed, $P_{j,i}\ne0$, and nondegeneracy of the
		trace pairing provides operators
		$X_1,\ldots,X_b\in M_{d_i}(E'_i)$ linking
		nonzero vectors in the images of the successive
		projectors so that
		$\operatorname{tr}(P_{1,i}X_1\cdots P_{b,i}X_b)\ne0$.
		The linear span of $\rho_i(S)$ is still the full
		matrix algebra after extension to $E'_i$.
		Multilinearity therefore prevents $F_i$ from
		vanishing on all of $S^b$.
		
		The functions $F_i$ are regular over their local
		coefficient fields. Lemma~\ref{lem:local-coefficients}
		and geometric irreducibility of $\mathcal Z^b$
		allow us to choose one list $h_1,\ldots,h_b\in S$
		such that
		\begin{equation}\label{eq:connectors}
			F_i(h_1,\ldots,h_b)\ne0\quad\text{for every }i.
		\end{equation}
		Set $X_n=g_1^{n!}h_1\cdots g_b^{n!}h_b\in S$.
		Lemma~\ref{lem:compact-return} gives, for every
		$i$ and every fixed $h\in\mathcal Z(\mathbb Q)$,
		\begin{equation}\label{eq:trace-limit}
			\left(\prod_j\beta_{j,i}\right)^{-n!}
			\operatorname{tr}\rho_i(X_nh)
			\longrightarrow
			F_i^h:={\operatorname{tr}}\bigl(P_{1,i}\rho_i(h_1)\cdots
			P_{b,i}\rho_i(h_b)\rho_i(h)\bigr).
		\end{equation}
		In particular,
		$F_i^{1}=F_i(h_1,\ldots,h_b)\ne0$.
		The set
		\[
		\Omega=\{h\in\mathcal Z:\ F_i^h\ne0\text{ for every }i\}
		\]
		is understood on rational points through these
		finitely many algebraic conditions, whose
		coefficients may lie in different local fields.
		Lemma~\ref{lem:local-coefficients} supplies a
		nonempty Zariski-open set defined over $\mathbb Q$ whose intersection
		with $S$ is contained in $S\cap\Omega$.
		In particular, $S\cap\Omega$ remains geometrically
		Zariski dense.
		
		For every geometric point $x\in\mathcal Z$, the
		open set $x^{-1}\mathcal U$ is nonempty and
		therefore meets $S\cap\Omega$. Hence the opens
		$\mathcal U h^{-1}$, for $h\in S\cap\Omega$,
		cover $\mathcal Z$. Noetherianity gives a finite
		subcover
		$\mathcal U a_1^{-1},\ldots,\mathcal U a_N^{-1}$,
		where $a_j\in S\cap\Omega$.
		For each $n$, choose $j(n)$ such that
		$X_na_{j(n)}\in\mathcal U$. On an infinite
		subsequence, $j(n)$ is constant, say equal to
		$j_0$. Since $F_i^{a_{j_0}}\ne0$ for all $i$,
		\eqref{eq:trace-limit} and
		\eqref{eq:beta-growth} yield, for all
		sufficiently large indices in this subsequence,
		\[
		|\operatorname{tr}\rho_i(X_na_{j_0})|_{E_i}
		=M_i^{n!}|F_i^{a_{j_0}}|_{E'_i}\longrightarrow\infty.
		\]
		Every $X_na_{j_0}$ is a positive word in the
		same semigroup $S$ and belongs to $\mathcal U$.
		No elements have been chosen independently in
		the different completions.
	\end{proof}
	
	We say that an element has \emph{double escape}
	if it has at least two eigenvalues of modulus
	greater than one, counted with multiplicity.
	In a unitary factor, reciprocity then gives at
	least two eigenvalues of modulus less than one.
	
	\begin{proposition}
		\label{prop:filtered-expansion}
		There are a subsemigroup $S\subset\Gamma_0$ that
		is $\mathbb Q$-Zariski dense in $\mathcal Z$ and
		a partition $W=A\sqcup B_u\sqcup B_c$ with the
		following properties:
		\begin{enumerate}[label=\textup{(\roman*)},leftmargin=*]
			\item Every element of $S$ is loxodromic at the
			distinguished place and has at least two
			expanding eigenvalues at each nonsplit real
			place where some $q_j$ is negative.
			\item Every element of $S$ has double escape in the
			factors of $A$, and none has double escape
			in $B_u\cup B_c$.
			\item The projection of $S$ is unbounded in
			$B_u$ and bounded in $B_c$.
			\item There is a single $t\in S$, all of whose
			positive powers are regular, such that
			$t_w\sim\operatorname{diag}(\alpha_w,R_w,\bar\alpha_w^{-1})$
			in $B_u$, with $|\alpha_w|>1$ and precompact
			powers of $R_w$. The powers of $t_w$ are
			precompact in $B_c$.
			\item There is a common strictly increasing sequence of positive integers
			$n_\ell$ with $R_w^{n_\ell}\to I$ in $B_u$
			and $t_w^{n_\ell}\to I$ in $B_c$.
		\end{enumerate}
		Once these data have been chosen, $S$ is kept
		fixed throughout the selection and return
		arguments.
	\end{proposition}
	\begin{proof}
		\emph{Filtering while preserving density.}
		Begin with $\Gamma_0$, viewed as a semigroup.
		In a real factor of rank at least $j$,
		Lemma~\ref{lem:proximal-semigroup} gives
		$g\in S$ whose $j$th exterior power is
		proximal. Algebraic regularity of this
		initial element is unnecessary: regularity
		of all powers of the final expanding element
		will be imposed later by a rational open
		set $\mathcal R$.
		If $\lambda$ is the dominant eigenvalue,
		then $\lambda^{-N}\rho(g^N)\longrightarrow
		P=u\ell$, with $\ell(u)=1$.
		For connecting elements $h_i$ satisfying
		$c_i=\ell(\rho(h_i)u)\ne0$, put
		$y=g^N$ and $z_i=g^Nh_ig^N$.
		Their exact normalizations satisfy
		\[
		\lambda^{-N}\rho(y)\longrightarrow P,\qquad
		(\lambda^{2N}c_i)^{-1}\rho(z_i)\longrightarrow P,
		\]
		because $P\rho(h_i)P=c_iP$.
		Thus all normalized matrices approach the
		same projector, with their scalar factors
		specified separately.
		Choose the $h_i$ so that the operators
		$\operatorname{Ad}(h_i)$ span the same
		rational vector space as the
		$\operatorname{Ad}(s)$ for elements of the
		current semigroup.
		Lemma~\ref{lem:local-coefficients} allows
		this choice within the indicated open set.
		
		For uniformity in a real factor, use an
		adapted basis in which a normalized matrix
		has blocks
		$\bigl(\begin{smallmatrix}a&b\\c&D\end{smallmatrix}\bigr)$
		with $|a-1|,\|b\|,\|c\|,\|D\|\le1/16$.
		On the ball $\|z\|\le1$, the projective map
		$F(z)=(c+Dz)/(a+bz)$ has norm at most $1/7$
		and derivative norm at most $4/49$.
		A product of length $r$ therefore contracts
		by at most $(4/49)^r$.
		Its fixed point gives a simple dominant
		eigenvalue, and the differential there
		identifies the ratios of the other
		eigenvalues to it. Every positive product
		is consequently proximal. Unitary
		reciprocity and proximality of the
		$j$th exterior power give at least $j$
		expanding eigenvalues when $j\le m/2$.
		
		The Zariski closure of a semigroup in
		an algebraic group is a group.
		Indeed, for an element $s$ of the
		semigroup, let
		$Z_s=\overline{\{s^n:n\ge1\}}^{\rm Zar}$.
		The descending chain of closed sets
		$Z_s\supset sZ_s\supset s^2Z_s\supset\cdots$
		stabilizes; translating gives $Z_s=sZ_s$.
		Since $s\in Z_s$, this implies $1\in Z_s$
		and then $s^{-1}\in Z_s$.
		The closure of the semigroup thus
		contains inverses and equals the closure
		of the group it generates.
		In particular, the Zariski closure
		$H'$ of the new semigroup is a group
		and contains $h_i=y^{-1}z_iy^{-1}$.
		The element $y$ has infinite order,
		so $\operatorname{Lie}(H')\ne0$.
		Stability under the $\operatorname{Ad}(h_i)$
		implies, by linearity, stability under
		all $\operatorname{Ad}(s)$ from the
		current semigroup. The stabilizer of
		this subspace is algebraic; density
		therefore makes $\operatorname{Lie}(H')$
		an ideal of $\operatorname{Lie}(\mathcal Z)$.
		Rational simplicity gives $H'=\mathcal Z$.
		Inverses are used only in this density
		argument and are never added to the
		positive generators.
		
		\emph{Imposing the real and finite-place spectral conditions.}
		Apply this construction first with $j=1$
		at the distinguished place and then
		with $j=2$ at each real place with a
		negative coefficient. The latter rank
		is available because the original form
		contains a hyperbolic plane, a line
		with coefficient $1$, and a negative
		line.
		At each finite factor where the current
		semigroup contains an element with
		double escape, repeat the construction
		with the $r$th exterior power, where
		$r$ is the number of expanding
		eigenvalues of such an element.
		Its dominant eigenvalue is unique
		and simple even if the element is
		not regular: it is the product of
		all expanding eigenvalues.
		The normalized powers converge to
		a rank-one projector. For sufficiently
		large $N$, all normalized matrices
		reduce to $\operatorname{diag}(1,0,\ldots,0)$
		in an adapted basis. Every product
		has the same reduction, since $P^2=P$,
		and Hensel's lemma gives proximality
		for every word.
		Lemma~\ref{lem:wedge-escape}, with
		$r\le\lfloor m/2\rfloor$ by unitary
		reciprocity, shows that each word
		retains at least $r$ expanding
		eigenvalues.
		The rational density argument is
		unchanged.
		After at most $|W|$ steps we obtain
		$A$ and $B$; properties imposed on
		every element are preserved on
		passing to subsemigroups.
		Now divide $B$ into $B_u$ and $B_c$.
		If $W$ is empty, this finite-place
		filtering is unnecessary and all
		three parts are empty.
		
		\emph{A common expanding element with regular powers.}
		To choose $t$ without replacing
		$S$, first fix the open set imposing
		regularity of all positive powers.
		The splitting field of an $m\times m$
		matrix over $K$ has absolute degree
		at most $D=m![K:\mathbb Q]$.
		The possible orders of roots of
		unity of such degree belong to
		the finite set
		$\{e:\varphi(e)\le D\}$.
		Let $M$ be a multiple of all
		these orders. The open set
		\[
		\mathcal R=\{g:\operatorname{disc}\chi_{g^M}\ne0\}
		\]
		is nonempty: in the original
		$\operatorname{SL}_m$ factor, one
		can choose torus eigenvalues that
		remain distinct after raising
		them to the $M$th power.
		The product of the Galois conjugates
		of the discriminant over
		$\mathbb Q$ defines a Zariski-open set
		defined over $\mathbb Q$ with the same condition
		on rational points. This product
		is nonzero by geometric integrality.
		If eigenvalues $\lambda_i,\lambda_j$ of $t$, counted with multiplicity
		and with $i\ne j$, satisfied $\lambda_i^n=\lambda_j^n$ for some $n>0$,
		then $\lambda_i/\lambda_j$ would be a root of unity
		whose order divides $M$. Hence
		$\lambda_i^M=\lambda_j^M$, contradicting the
		defining discriminant condition for $\mathcal R$.
		Thus membership in $\mathcal R$
		forces every positive power
		to be regular semisimple.
		
		The natural representations in
		$B_u$ are absolutely irreducible
		by Lemma~\ref{lem:source-models},
		and the image of the same $S$
		is unbounded by definition.
		Apply Lemma~\ref{lem:common-trace}
		to these representations and
		to $\mathcal R$. It gives
		$t\in S\cap\mathcal R$ with
		$|\operatorname{tr}t_w|_w>1$
		simultaneously for every
		$w\in B_u$.
		The connecting elements in that lemma lie in the same $S$;
		no topological density is needed.
		
		In $B_u$, absence of double
		escape and the trace inequality
		give exactly one expanding and
		one contracting eigenvalue.
		They are simple and belong to
		$E_w$ by uniqueness of the expanding
		and contracting eigenvalues and
		Galois invariance. Their eigenlines
		are isotropic and pair
		nontrivially; their orthogonal
		complement $W_w$ has dimension
		$m-2$. Regularity makes $R_w$
		semisimple, and its eigenvalues
		are units, so its powers are
		precompact. In $B_c$, boundedness
		of $S$ suffices.
		If $B_u$ is empty, choose
		$t\in S\cap\mathcal R$ directly.
		\emph{A common return sequence.}
		The closures of the powers
		of the $R_w$ and $t_w$ in
		these factors are compact
		groups: recurrence in a
		compact closure supplies
		the inverses.
		Their finite product is
		compact and totally disconnected.
		The factorial powers of its
		diagonal element converge to
		the identity, since they
		eventually become trivial
		in each finite quotient
		associated with a basis
		of open subgroups.
		We may take $n_\ell=\ell!$
		for $\ell\ge2$.
		If there are no residual
		factors, the condition is
		vacuous and $n_\ell=\ell$
		suffices.
	\end{proof}
	
	\subsection{Exterior-power return after the connectors have been fixed}
	
	The next proposition relates products in the transverse unitary groups
	to spectra of words in the original group. The return sequence is fixed
	first, the connecting elements are chosen next, and the exponent is
	chosen last. Only the final threshold depends on the connecting elements.
	
	\begin{proposition}
		\label{prop:return}
		Let $E/F$ be a quadratic extension of nonarchimedean local fields
		of characteristic zero, with involution $\tau$. Relative to
		$Ee_+\oplus W\oplus Ee_-$, let
		\[
		h=\begin{pmatrix}0&0&1\\0&H_W&0\\1&0&0\end{pmatrix},\quad
		t=\operatorname{diag}(\alpha,R,\tau(\alpha)^{-1})\in\mathrm{SU}(h)(F),
		\]
		where $|\alpha|>1$ and the powers of $R$ are precompact. First fix
		a sequence of positive integers $n_\ell\to\infty$ such that $R^{n_\ell}\to I$.
		Then fix any nonempty finite list $h_1,\ldots,h_b\in\mathrm{SU}(h)(F)$, $b\ge1$,
		with $a_j=(h_j)_{++}\ne0$, and set
		\[
		C_j=(h_j)_{WW}-(h_j)_{W+}a_j^{-1}(h_j)_{+W},\quad
		C=C_b\cdots C_1,\quad A=\prod_j a_j,\quad
		\gamma_\ell=t^{n_\ell}h_b\cdots t^{n_\ell}h_1.
		\]
		Then $C_j\in\mathrm U(H_W)(F)$. For sufficiently large $\ell$,
		$\gamma_\ell$ has simple roots
		$\lambda_{+,\ell},\lambda_{-,\ell}\in E$ satisfying
		\[
		\alpha^{-bn_\ell}\lambda_{+,\ell}\to A,\qquad
		\lambda_{-,\ell}=\tau(\lambda_{+,\ell})^{-1}\to0.
		\]
		The monic quotient
		\[
		Q_\ell(X)=\frac{\chi_{\gamma_\ell}(X)}
		{(X-\lambda_{+,\ell})(X-\lambda_{-,\ell})}
		\]
		converges coefficientwise to $\chi_C(X)$. The threshold depends
		on the fixed list; no uniformity over all list lengths is asserted.
		For finitely many places, a common return sequence can be chosen
		before the connectors are selected.
	\end{proposition}
	
	\begin{proof}
		The first column of $h_j$, divided by $a_j$, is
		$(1,v,d)^t$, where $d+\tau(d)+v^*H_Wv=0$. The matrix
		\[
		n(v,u)=\begin{pmatrix}
			1&0&0\\v&I&0\\u-\tfrac12v^*H_Wv&-v^*H_W&1
		\end{pmatrix},\qquad u=d+\tfrac12v^*H_Wv,
		\]
		is unitary, and $u+\tau(u)=0$. The matrix $n(v,u)^{-1}h_j$
		stabilizes $Ee_+$, and its central block is precisely $C_j$.
		Its action on the Hermitian quotient $e_+^\perp/Ee_+$ shows that
		$C_j$ is unitary.
		
		For $1\le r\le m-1$, let $P_r$ denote the projection of
		$\bigwedge^r(E^m)$ onto $e_+\wedge\bigwedge^{r-1}W$.
		The exterior-power decomposition gives
		\[
		\alpha^{-n_\ell}\bigwedge^r t^{n_\ell}\longrightarrow P_r.
		\]
		The terms that do not contain $e_+$, or that also contain $e_-$,
		acquire a factor tending to zero, while the powers of $R$ remain
		bounded. The formula for minors with pivot $a_j$ shows that the
		restriction of $P_r(\bigwedge^rh_j)P_r$ to its image is
		$a_j\bigwedge^{r-1}C_j$. If $e_r$ denotes the $r$th elementary
		symmetric function of the roots, multiplication and the trace give
		\begin{equation}
			\alpha^{-bn_\ell}e_r(\gamma_\ell)\longrightarrow A e_{r-1}(C),
			\qquad
			\alpha^{-bn_\ell}\chi_{\gamma_\ell}(X)
			\longrightarrow -AX\chi_C(X).
			\label{eq:exterior-return}
		\end{equation}
		Put $c_\ell=\alpha^{-bn_\ell}$. The monic polynomial
		$c_\ell^m\chi_{\gamma_\ell}(c_\ell^{-1}Y)$ converges to
		$Y^{m-1}(Y-A)$. Hensel's lemma at the simple root $A\ne0$
		gives $Y_\ell\in E$ with $Y_\ell\to A$, and hence the root
		$\lambda_{+,\ell}=c_\ell^{-1}Y_\ell$.
		
		To control the remaining roots without assuming that they are
		bounded, first divide by this exact factor:
		\[
		c_\ell\chi_{\gamma_\ell}(X)
		=(c_\ell X-Y_\ell)P_\ell(X),\quad
		P_\ell=\chi_{\gamma_\ell}/(X-\lambda_{+,\ell}).
		\]
		Comparison of coefficients, starting with the constant term and
		using $Y_\ell\to A\ne0$, gives $P_\ell\to X\chi_C$.
		Since $\det C\ne0$, the limit has a simple root at zero.
		Hensel's lemma gives a unique small root of $P_\ell$. Unitarity
		forces it to be $\tau(\lambda_{+,\ell})^{-1}$. Continuous division
		by this second monic factor gives $Q_\ell\to\chi_C$.
		
		The compact-group return argument in the proof of
		Proposition~\ref{prop:filtered-expansion}, applied to the finite
		product of precompact blocks, gives a common sequence.
		The connectors are selected afterwards; only the final threshold
		depends on them.
	\end{proof}
	
	\subsection{The compact closure and ambient open sets}
	
	The return formula leaves us to choose the connecting elements.
	We must approximate the required spectral targets by products of their
	transverse blocks in $B_u$ and full matrices in $B_c$. We therefore
	study the closure of these products near the identity.
	Consider source data as in \eqref{eq:source-product}: a
	rationally Zariski-dense semigroup $S\subset Z^\circ(\mathbb Q)$,
	a finite partition $A\sqcup B_u\sqcup B_c$, and an element $t\in S$.
	Assume that no element of $S$ has two expanding roots at a place in
	$B_u$, and that $t$ has the form in
	Proposition~\ref{prop:return}; at the places in $B_c$, the image
	of $S$ is precompact. Factors corresponding to extensions of the
	same place $w$ are identified with a single source. Let $\Omega$ be
	the cell on which all pivots $(h_w)_{++}$ are nonzero, and set
	\[
	\Phi(h)=((C_w(h))_{w\in B_u},(h_w)_{w\in B_c}),\quad
	\mathcal T=\langle\Phi(S\cap\Omega)\rangle_+,\quad
	C=\overline{\mathcal T}.
	\]
	Here $\langle\,\cdot\,\rangle_+$ denotes products of positive length.
	
	\begin{proposition}
		\label{prop:effective-compact}
		Under the preceding assumptions, $C$ is a compact group. Write
		\[
		J_p=\prod_{\substack{w\in B_u\\w\mid p}}
		\operatorname{Res}_{f_w/\mathbb Q_p}\mathrm U(W_w)
		\ \times\!
		\prod_{\substack{w\in B_c\\w\mid p}}
		\operatorname{Res}_{f_w/\mathbb Q_p}\mathrm{SU}(h_w).
		\]
		The projection $C_p$ is jointly Zariski-dense in $J_p$. There are
		open identity neighborhoods
		$V_p\subset J_p^{\mathrm{der}}(\mathbb Q_p)$ such that
		$\prod_pV_p\subset C$.
	\end{proposition}
	
	\begin{proof}
		Fix an arbitrary product in $\mathcal T$ and its connectors.
		By Proposition~\ref{prop:return}, the characteristic polynomials
		of its transverse components are limits of the central characteristic
		polynomials of words in $S$, after separating one expanding root.
		Since no element of $S$ has two expanding roots,
		every approximating central polynomial has all its roots of absolute
		value one. A monic polynomial whose determinant is a unit has all
		its roots of absolute value one if and only if its coefficients are
		integral: one implication follows from the ultrametric inequality,
		and the other from integrality of the roots and the fact that their
		product is a unit. This set of polynomials is closed. Thus every
		transverse product in $\mathcal T$ has all its roots of absolute
		value one.
		
		The joint map $\Phi$ is dominant onto $J_p$. In a transverse
		factor, for $B\in\mathrm U(W)$ with $d=\det B\ne-1$, the section
		\[
		B\longmapsto\operatorname{diag}(a,B,\tau(a)^{-1}),\qquad
		a=(1+d^{-1})/2,
		\]
		has determinant one because $a/\tau(a)=d^{-1}$, and its Schur
		complement is $B$. In $B_c$, use the identity map. The product
		of sections is taken over the distinct source factors in
		\eqref{eq:source-product}; no independent sections are assigned
		to correlated extensions of a place.
		Lemma~\ref{lem:local-coefficients} gives the density of
		$\Phi(S\cap\Omega)$ in this product.
		
		The natural representation of each $\mathrm U(W_w)$ is absolutely
		irreducible. By Burnside's theorem, the linear span of the projection
		of $\mathcal T$ is the full matrix algebra. Choose a basis
		$T_1,\ldots,T_{r^2}$ from this projection. For $T\in\mathcal T$,
		every $\operatorname{tr}(TT_i)$ is integral, since $TT_i$ is another
		product with all roots of absolute value one. The matrix
		$(\operatorname{tr}(T_iT_j))_{i,j}$ is invertible by nondegeneracy
		of the trace pairing. Solving this linear system uniformly bounds
		the entries of $T$. Since its determinant is a unit, its inverse
		is bounded as well. The factors in $B_c$ are already bounded.
		Thus $C$ is compact in a linear group, and recurrence of powers
		makes it a group.
		
		Let $D_p=\overline{[C_p,C_p]}$. Its Zariski closure is
		$J_p^{\mathrm{der}}$, by density and the commutator map. Its Lie
		algebra is stable under $C_p$, and the stabilizer of a subspace is
		algebraic. This Lie algebra is therefore an ideal of
		$\operatorname{Lie}(J_p^{\mathrm{der}})$. Each restriction of scalars
		of an absolutely simple local Lie algebra is simple over $\Q_p$:
		after scalar extension its simple ideals are permuted transitively
		by Galois, so a rational ideal contains either all or none of them.
		Write the resulting decomposition as
		$\mathfrak j=\bigoplus_a\mathfrak j_a$. An ideal of $\mathfrak j$
		is a sum of these simple factors. If $\Lie(D_p)$ omitted one factor,
		the image of $D_p$ in the corresponding adjoint factor would have
		zero-dimensional Lie algebra. This image is compact and closed;
		the closed subgroup theorem makes it a discrete $p$-adic Lie group,
		and compactness then makes it finite. Its Zariski density in the
		positive-dimensional adjoint factor rules this out. Hence every
		factor occurs and
		$\operatorname{Lie}D_p=\operatorname{Lie}J_p^{\mathrm{der}}$.
		The $p$-adic closed subgroup theorem and an exponential chart
		then give an identity neighborhood in the derived group contained
		in $D_p$; see \cite[\S2, Propositions 2.2--2.3]{Glockner2016}.
		
		To separate distinct primes, choose open normal pro-$p$ subgroups
		$U_p\subset C_p$, using \cite[Theorem 3.1]{Glockner2016} and
		passing to the normal core, and put $C'=C\cap\prod_pU_p$.
		Its projections $H_p$ have finite index in $C_p$ and are pro-$p$.
		A common continuous quotient of a pro-$p$ group and a pro-$q$
		group, with $p\ne q$, is trivial: each of its finite quotients
		would have order both a power of $p$ and a power of $q$.
		At each inductive step, a finite quotient of the product of the
		other prime factors has order prime to $p$. Profinite Goursat
		therefore gives $C'=\prod_pH_p$.
		Each $H_p$ remains Zariski-dense because $J_p$ is connected.
		Applying the preceding Lie algebra argument to $H_p$ gives the
		neighborhoods $V_p$ with $\prod_pV_p\subset C$.
	\end{proof}
	
	\begin{corollary}
		\label{cor:one-word}
		For each residual factor, choose a target in an identity neighborhood
		in the derived group contained in $C$, together with an ambient open
		neighborhood preserving the required polynomial or spectral type.
		There is a single product $\Phi(h_b)\cdots\Phi(h_1)$, with $b\ge1$, in these
		neighborhoods. Once this list has been fixed, a sufficiently large
		index in the return sequence yields one word
		$t^{n_\ell}h_b\cdots t^{n_\ell}h_1\in S$ with all the prescribed
		properties: in $B_u$ through its central factor, and in $B_c$
		through its full characteristic polynomial.
	\end{corollary}
	\begin{proof}
		The tuple of targets belongs to $C$. The product of their ambient
		open neighborhoods meets $C$ in a nonempty relatively open set.
		Density of $\mathcal T$ supplies a nonempty finite product in that set;
		this product is not required to lie in the derived group.
		Once its letters have been fixed, the convergence in
		Proposition~\ref{prop:return}, the return of $t$ in the bounded
		factors, and finiteness of the set of places give a single threshold.
		The argument uses the product of the values of $\Phi$, never the
		false identity $\Phi(h_b\cdots h_1)=\Phi(h_b)\cdots\Phi(h_1)$.
	\end{proof}
	
	\subsection{Norms from a word with the required local properties}
	
	For a regular semisimple loxodromic word $\gamma$, we seek local
	spectral conditions that make every $q_j$ a norm from its splitting
	field $L=K(E(\gamma))$ to $L^c$. The signed Galois action keeps track
	of actual complex conjugation even when it is not central in
	$\operatorname{Gal}(L/k)$.
	
	\begin{lemma}
		\label{lem:signed-galois}
		Under the hypotheses of Theorem~\ref{thm:norm-numeric}, let
		$\gamma\in\mathrm{SU}(H)(k)$ be regular semisimple and loxodromic,
		with distinct roots $\lambda_1,\ldots,\lambda_m$. Then
		$L=K(\lambda_1,\ldots,\lambda_m)$ is Galois over $k$.
		For $\sigma\in\operatorname{Gal}(L/k)$, put
		$\varepsilon_\sigma=1$ when $\sigma|_K=1$ and
		$\varepsilon_\sigma=-1$ otherwise. There is a unique permutation
		$\pi_\sigma$ such that
		\[
		\sigma(\lambda_i)=\lambda_{\pi_\sigma(i)}^{\varepsilon_\sigma}.
		\]
		The map $\sigma\mapsto(\pi_\sigma,\varepsilon_\sigma)$ is an
		injective homomorphism into $S_m\times\{1,-1\}$. For actual
		complex conjugation, $\pi_c$ is a single transposition; the same
		holds for every conjugate of $c$.
	\end{lemma}
	\begin{proof}
		If $p(X)=\det(X-\gamma)$, unitarity and determinant one give
		$\overline p(X)=(-1)^mX^mp(X^{-1})$. Every $k$-embedding of the
		splitting field sends the roots to roots of $p$ or their inverses,
		according to its action on $K$. It therefore preserves $L$, which
		is normal and separable. Distinctness of the roots gives uniqueness
		of $\pi_\sigma$; composition of the identities shows that the
		permutations compose and the signs multiply. The kernel fixes
		$K$ and every root, so it is trivial. In signature $(m-1,1)$,
		there is a unique pair of roots off the unit circle,
		$\lambda_- =c(\lambda_+)^{-1}$; all other roots satisfy
		$c(\lambda)=\lambda^{-1}$. The signed action of $c$ interchanges
		this pair and fixes the remaining labels. Cycle type is preserved
		by conjugation in $S_m\times\{1,-1\}$, without any assumption
		that $c$ is central.
	\end{proof}
	
	\begin{lemma}
		\label{lem:even-local-norm}
		Let $K/k$ be a quadratic extension of number fields, let $L/k$
		be a finite Galois extension containing $K$, let $c\in\operatorname{Gal}(L/k)$ be an
		involution acting nontrivially on $K$, and put $M=L^c$.
		If $v$ is a finite place of $k$, $K_v=K\otimes_k k_v$ is a field, and
		$[L_{\mathfrak w}:K_v]$ is even for one extension of $v$, then
		every $q\in k_v^\times$ is a norm from
		$L\otimes_M M_u/M_u$ for every place $u\mid v$ of $M$.
	\end{lemma}
	\begin{proof}
		The degree is the same at all extensions of $v$, since $L/k$
		is Galois. If $c$ interchanges the places above $u$, the local
		algebra is $M_u\times M_u$ and its norm is surjective. If $c$
		fixes $\mathfrak w$, it acts on $L_{\mathfrak w}$ with fixed
		field $M_u$, and
		\[
		L_{\mathfrak w}=K_vM_u,\quad K_v\cap M_u=k_v,\quad
		[M_u:k_v]=[L_{\mathfrak w}:K_v].
		\]
		The cyclic algebra class $(K_v/k_v,q)$ has invariant $0$ or
		$1/2$. Restriction to $M_u$ multiplies this invariant by the
		even degree above \cite[III, equation (29)]{MilneCFT2020}, so the
		restricted class vanishes. The splitting criterion for the cyclic
		algebra says precisely that $q$ is a norm. An auxiliary quadratic
		extension contained in $L_{\mathfrak w}$ can be used to establish
		evenness; it need not be stable under $c$.
	\end{proof}
	
	\begin{proposition}
		\label{prop:norm-assembly}
		Retain the hypotheses of Theorem~\ref{thm:norm-numeric}, and let
		$\gamma\in\Gamma$ be regular semisimple and loxodromic. Let
		$P_f$ be the union of the finite places at which at least one
		$q_j$ is not a norm from $K/k$. Assume that this same word
		satisfies the following conditions.
		\begin{enumerate}[label=\textup{(\roman*)}]
			\item At every $v\in P_f$, one of the following holds:
			there are at least two expanding and two contracting roots;
			or, if $m=4$, the degree $[L_{\mathfrak w}:K_v]$ is even;
			or, if $m\ge5$, the roots split in $K_v$ and the local
			involution has signed permutation type equal to the identity
			or to two disjoint transpositions.
			\item At every real place $v$ with $K_v=\mathbb C$ and
			$q_j(v)<0$ for some $j$, there are at least two expanding
			and two contracting roots.
		\end{enumerate}
		Then every $q_j$ is a norm from $L/L^c$, where
		$L=K(E(\gamma))$ is the splitting field of this one word.
	\end{proposition}
	\begin{proof}
		By Lemma~\ref{lem:signed-galois}, $L/k$ is Galois and
		$L/M$, with $M=L^c$, is quadratic, with $L=KM$ and
		$K\cap M=k$. If $c$ fixes a place of $L$, it preserves the
		corresponding absolute value. The $m-2$ identities
		$c(\lambda)=\lambda^{-1}$ for the distinguished central roots
		force those roots to have absolute value one there. At most
		two roots can therefore lie off the unit circle. Double escape
		thus rules out every place fixed by $c$ above the corresponding
		$v$; the local extension of $L/M$ is split, and its norm is
		surjective. At an obstructing real place, the same contradiction
		prevents $M_u$ from being real: all places above it are complex,
		and the local quadratic algebra is split.
		
		For the even-degree alternative, apply
		Lemma~\ref{lem:even-local-norm}. In the signed-permutation
		alternative, the roots lie in $K_v$, so
		$L_{\mathfrak w}=K_v$. Its decomposition group over $k_v$
		has order two, and its nontrivial element has permutation type
		equal to the identity or to a double transposition. Neither is
		conjugate to $c$, whose permutation type is a single
		transposition. The decomposition groups of all places above
		$v$ are conjugate. Thus no such place is fixed by $c$, and
		$L/M$ is again split locally.
		
		At every remaining place, each $q_j$ was already a norm from
		the algebra $K\otimes_k k_v$. Base change through
		$K\otimes_k M\simeq L$ carries these norm witnesses to every
		completion of $M$. This covers every place, for every
		coefficient and the same quadratic extension of number fields.
		The Hasse norm theorem
		\cite[VIII, Theorem 3.1]{MilneCFT2020} completes the proof.
	\end{proof}
	
	\subsection{Targets and correlations for simultaneous norm realization}
	
	\begin{lemma}
		\label{lem:original-parity}
		Let $E/F$ be a quadratic extension of nonarchimedean local fields
		of characteristic zero, let $k'/F$ be finite, and suppose that
		$K'=Ek'$ is a quadratic field extension of $k'$. If
		$H'=\mathbb H\perp\langle1,q\rangle$ is similar to a
		four-dimensional Hermitian form $h$ over $E/F$ after base change,
		conjugation, or duality, then
		\[
		[(K'/k',q)]=\operatorname{Res}_{k'/F}[(E/F,-\det h)].
		\]
		If $q$ is not a norm, then $[k':F]$ is odd.
	\end{lemma}
	\begin{proof}
		The determinant of $H'$ is $-q$. A change of basis multiplies
		the determinant by a norm, while a similarity multiplies it by
		$a^4=N(a^2)$. Conjugation leaves the determinant in $F$
		unchanged, and duality inverts it, changing its class by the
		square $(\det h)^2$, which is also a norm. This gives equality
		of the classes modulo norms, and hence of the cyclic algebra
		classes. If $q$ is not a norm, restriction sends an invariant
		equal to $0$ or $1/2$ to $1/2$, forcing the degree to be odd.
		The argument uses the coefficient of the original form; it
		does not apply to an independently chosen local coefficient $q$.
	\end{proof}
	
	\begin{lemma}
		\label{lem:binary-target}
		Let $E/F$ be a quadratic extension of nonarchimedean local fields
		of characteristic zero, let $R=F(\sqrt d)\ne E$ be quadratic,
		and put $M=ER$. For every binary Hermitian form $h_2$, every
		identity neighborhood in $\mathrm{SU}(h_2)(F)$ contains a
		regular element with splitting field $M/E$. This property
		persists in an ambient open neighborhood in
		$\mathrm U(h_2)(F)$. For a four-dimensional form, the same
		splitting field is obtained by adjoining two distinct roots in $E$.
	\end{lemma}
	\begin{proof}
		Up to similarity, diagonalize $h_2$ as $\operatorname{diag}(1,r)$.
		Its special unitary group is the reduced norm-one group of the
		quaternion algebra
		\[
		D=\left\{\begin{pmatrix}u&-r\tau(v)\\v&\tau(u)\end{pmatrix}:u,v\in E\right\}.
		\]
		The adjoint involution of $h_2$ is the canonical involution of
		$D$. Every quadratic extension $R/F$ embeds in $D$: if $D$
		is split, use the regular representation; if it is a division
		algebra, the degree-two extension kills its invariant $1/2$,
		and the maximal-subfield criterion gives the embedding
		\cite[IV, Corollary 3.7 and Remark 4.4(c)]{MilneCFT2020}.
		The image $X$ of $\sqrt d$ satisfies $X^2=dI$, has trace
		zero, and satisfies $X^\dagger=-X$. For sufficiently small
		nonzero $s\in F$, its Cayley transform is
		\[
		B_s=(I+sX)(I-sX)^{-1}
		=\frac{(1+ds^2)I+2sX}{1-ds^2},\qquad
		\operatorname{disc}\chi_{B_s}
		=d\left(\frac{4s}{1-ds^2}\right)^2.
		\]
		It is special unitary and tends to the identity. Its roots
		are $(1+s\sqrt d)/(1-s\sqrt d)$ and its inverse, and they
		generate $M/E$. The square class of a nonzero discriminant
		is locally constant because $E^{\times2}$ is open, including
		in residue characteristic two. Hensel's lemma preserves the
		quadratic factor in a coefficient neighborhood; its preimage
		is taken in $\mathrm U(h_2)$, with no determinant constraint.
		In dimension four, separate a second diagonal plane and use
		norm-one roots $\nu,\nu^{-1}\in E$, close to one and distinct
		from the first two roots. The coprime factorization persists
		by Hensel's lemma. The parameter can be chosen in sufficiently
		deep principal units, excluding torsion quotients of roots.
	\end{proof}
	
	\paragraph{One quadratic extension for each source.}
	By Lemma~\ref{lem:original-parity}, every obstructing original
	extension $v\mid w$ has odd degree over $F_w$. Choose a single
	quadratic extension $R_w^{\mathrm{aux}}/F_w$ distinct from
	$E_w$, and put $M_w=E_wR_w^{\mathrm{aux}}$. Then
	$M_w\cap K_v=E_w$ for every such extension, and
	$[M_wK_v:K_v]=2$. The binary target in $B_u$, or the
	four-dimensional target in $B_c$, is chosen once for each $w$.
	Corollary~\ref{cor:one-word} realizes these neighborhoods
	with a single word. Transport by conjugation or duality
	preserves the field, as does multiplication by a nonzero
	scalar in the base field. This gives the even degree required
	by Proposition~\ref{prop:norm-assembly}.
	
	\begin{lemma}
		\label{lem:signed-target}
		Let $E/F$ be a quadratic extension of nonarchimedean local fields
		of characteristic zero, with nontrivial involution $\tau$, and let $W$ be a Hermitian space of
		dimension $r\ge3$. Every identity neighborhood in
		$\mathrm{SU}(W)(F)$ contains a regular element with roots in
		$E$ for which the type
		$\tau(\lambda_i)=\lambda_{\pi(i)}^{-1}$ is a transposition.
		In every dimension $r\ge2$, there are targets of identity
		type. Both types persist in ambient open neighborhoods in
		$\mathrm U(W)(F)$.
	\end{lemma}
	\begin{proof}
		The quadratic form $v\mapsto h(v,v)$ over $F$ is nondegenerate
		of dimension $2r\ge6$. It is isotropic by
		\cite[VIII, Corollary 3.10]{MilneCFT2020}, including at dyadic
		places. Split off $W=\mathbb H\perp W'$, and use roots
		$a,a^{-1}\in F^\times$ on $\mathbb H$, distinct and close
		to one. Their signed type is a transposition. On a
		diagonalization of $W'$, choose roots $\exp(b_jx)$, where
		$\tau(x)=-x$, $x$ is nonzero and sufficiently small, and the $b_j\in\mathbb Q$
		are distinct and sum to zero. If $\dim W'=1$, use the
		root one. These roots have norm one and product one.
		They do not coincide with $a,a^{-1}$, since a norm-one
		element of $F$ is $\pm1$. For identity type, use this
		diagonal construction on all of $W$.
		
		Choose disjoint multiplicative balls around the roots that
		are permuted by involution followed by inversion according
		to $\pi$. Hensel's lemma gives a unique simple root in
		each ball for every sufficiently close polynomial.
		Unitary reciprocity forces the same permutation. The
		preimage of the coefficient neighborhood in $\mathrm U(W)$
		is an ambient open set. Under natural or dual transport,
		or multiplication by a norm-one scalar $z$, the roots
		$\mu_i=z\iota(\lambda_i)^{\varepsilon}$ satisfy
		$\tau'(\mu_i)=\mu_{\pi(i)}^{-1}$, so the type is preserved.
	\end{proof}
	
	\paragraph{Residual types in higher dimension.}
	For $m\ge5$, the transverse space in $B_u$ has dimension
	$m-2\ge3$. A central target of transposition type, together
	with the extreme pair supplied by return, gives two disjoint
	transpositions. In $B_c$, use identity type. These are the
	residual conditions of Proposition~\ref{prop:norm-assembly}.
	Proposition~\ref{prop:filtered-expansion} already supplies
	loxodromy at the distinguished place and double escape in $A$
	and at the obstructing real places. Thus the targets complete
	the local spectral conditions required for norm realization.
	We now make the simultaneous choice in the prescribed group.
	
	\begin{proof}[Proof of Theorem~\ref{thm:norm-numeric}]
		Fix $K/k$, $H$, $\Gamma$, and all the coefficients $q_j$.
		Apply Proposition~\ref{prop:rational-source} and
		Lemma~\ref{lem:source-models}, followed by
		Proposition~\ref{prop:filtered-expansion}, to obtain the source
		$\mathcal Z$, a single semigroup $S\subset\Gamma_0$, the partition
		$A\sqcup B_u\sqcup B_c$, an expanding element $t$, and a common
		return sequence. Every element of $S$ is loxodromic at the distinguished
		place and has double escape in $A$ and at the obstructing real places.
		Keep $S$ fixed throughout the selection.
		
		If $B_u\cup B_c$ is empty, take $\gamma=t$. All its powers are
		regular, and every obstructing finite place belongs to $A$.
		Proposition~\ref{prop:norm-assembly} gives the norms, including when
		there are no obstructing finite places.
		
		Otherwise, construct $\Phi$, $\mathcal T$, and $C$ from these data.
		Proposition~\ref{prop:effective-compact} gives a product of identity
		neighborhoods in the derived factors contained in $C$.
		For $m=4$, choose in these neighborhoods the binary targets in $B_u$
		and four-dimensional targets in $B_c$ of Lemma~\ref{lem:binary-target},
		using one auxiliary quadratic extension at each source as above.
		Lemmas~\ref{lem:source-models} and~\ref{lem:original-parity} give
		odd degree at every obstructing original extension, so the quadratic
		spectral extension remains quadratic after all these base changes.
		For $m\ge5$, Lemma~\ref{lem:signed-target} gives central targets
		of transposition type in $B_u$ and identity type in $B_c$.
		In both cases, the spectral conditions persist in ambient open
		neighborhoods; the approximating product need not lie in the derived group.
		
		Corollary~\ref{cor:one-word} selects a single list
		$h_1,\ldots,h_b\in S\cap\Omega$, followed by a common return
		exponent, such that
		\[
		\gamma=t^{n_\ell}h_b\cdots t^{n_\ell}h_1\in S\subset\Gamma.
		\]
		The list is fixed before the threshold and realizes all targets
		simultaneously. For $m=4$, the quadratic factor gives even local
		spectral degree at every obstructing original extension.
		For $m\ge5$, all local roots split: the type in $B_c$ is the
		identity, while the extreme pair in $B_u$ adds a disjoint transposition
		to the central one. Natural or dual transport preserves these conditions.
		
		The target polynomials are separable. In $B_u$, their central unit
		roots are separated from the two extreme roots; in $B_c$, the full
		polynomial is separable. Regularity at one residual factor therefore
		makes the discriminant of the original matrix nonzero, so $\gamma$
		is regular semisimple. Its loxodromy and the double-escape conditions
		hold because $\gamma\in S$. Proposition~\ref{prop:norm-assembly}
		now applies to this same word and gives \eqref{eq:simultaneous-norms}.
	\end{proof}
	
	\subsection{Powers of a norm-realizing element and persistence of the splitting field}
	\label{subsec:powers}
	
	Regularity of every power of the auxiliary expanding element
	$t$ does not by itself imply the same property for
	$t^{n_\ell}h_b\cdots t^{n_\ell}h_1$. The following choice
	of targets ensures this property for the final word by
	means of a local open condition.
	
	\begin{lemma}
		\label{lem:deep-units}
		Let $M/\mathbb Q_p$ be finite, let $v_M$ be its valuation
		with value group $\mathbb Z$, and put $e_M=v_M(p)$.
		If the integer $a$ satisfies $a>e_M/(p-1)$, then the group
		\[
		U_M(a)=1+\mathfrak p_M^a
		\]
		is torsion-free. More precisely, for
		$u\in U_M(a)\setminus\{1\}$ and every integer $n\ge1$,
		\begin{equation}
			v_M(u^n-1)=v_M(u-1)+e_Mv_p(n)<\infty.
			\label{eq:deep-units}
		\end{equation}
		In particular, two distinct elements of $U_M(a)$ have
		distinct $n$th powers for every $n\ge1$.
	\end{lemma}
	\begin{proof}
		Write $u=1+x$, with $s=v_M(x)\ge a$. If $p\nmid b$,
		the term $bx$ in $(1+x)^b-1$ has valuation $s$, whereas
		all remaining terms have valuation at least $2s>s$.
		Thus $v_M(u^b-1)=s$. In $(1+x)^p-1$, the term $px$
		has valuation $e_M+s$. The intermediate terms have
		valuation at least $e_M+2s$, since $p$ divides
		$\binom pj$ for $1\le j<p$, while the last term has
		valuation $ps>e_M+s$. The minimum is unique, including
		when $p=2$, and $v_M(u^p-1)=e_M+s$. Iteration for
		$n=p^rb$ proves \eqref{eq:deep-units}. The quotient
		of two elements of $U_M(a)$ again lies in this group,
		giving the final assertion. This elementary proof agrees
		with the torsion-freeness obtained from the injective local
		logarithm \cite[Theorems 8.5, 8.7, and 8.13]{ConradPadic}.
	\end{proof}
	
	\begin{lemma}
		\label{lem:deep-targets}
		The targets in Lemmas~\ref{lem:binary-target}
		and~\ref{lem:signed-target} can be chosen within
		any prescribed identity neighborhood so that all their
		roots are distinct and belong to one group $U_M(a)$
		from Lemma~\ref{lem:deep-units}. For the binary
		target and its four-dimensional extension, take $M=ER$;
		for the signed targets, take $M=E$. There is an ambient
		open neighborhood of the target in the corresponding
		unitary group that preserves both its original spectral
		condition and the absence of torsion quotients between
		distinct roots.
	\end{lemma}
	\begin{proof}
		The parameters in Lemmas~\ref{lem:binary-target}
		and~\ref{lem:signed-target} can be made arbitrarily small, placing
		all roots in $U_M(a)$ while keeping the target in the prescribed
		identity neighborhood. The binary Cayley roots remain distinct for
		small nonzero parameters. In dimension four, choose the additional
		pair $\nu,\nu^{-1}\in E\cap U_M(a)$ avoiding the finitely many
		collisions. Infinitely many norm-one choices are supplied by
		$(1+tz)/(1-tz)$, where $0\ne z$, $\tau(z)=-z$, and $t\in F$
		is sufficiently small.
		
		For the signed targets, retain the zero-sum parameters that ensure
		determinant one. A one-dimensional complement contributes the root one;
		the hyperbolic roots $a_0,a_0^{-1}\in F$ avoid $\pm1$ and therefore
		cannot coincide with norm-one roots. Identity-type targets use the
		same diagonal construction.
		
		Choose disjoint balls around the roots inside $U_M(a)$. Hensel's
		lemma over $M$ gives a coefficient neighborhood over $E$ with one
		simple root in each ball. The roots remain distinct in $M$ and their
		quotients lie in the torsion-free group $U_M(a)$.
		Intersect with the neighborhood preserving the original quadratic
		discriminant or signed type. Its preimage under the characteristic-polynomial
		map is the required ambient open set in $\mathrm U(W)(F)$, or in
		the full group for a bounded factor; the approximating product need
		not lie in the derived group.
	\end{proof}
	
	\begin{lemma}
		\label{lem:power-field}
		Let $K$ be a field of characteristic zero, and let
		$g\in\mathrm{GL}_m(K)$ have distinct roots
		$\lambda_1,\ldots,\lambda_m$. If no quotient
		$\lambda_i/\lambda_j$, with $i\ne j$, is a root of
		unity, then, for every integer $n\ge1$, the matrix
		$g^n$ is regular semisimple and
		\begin{equation}
			K(E(g^n))=K(E(g)).
			\label{eq:power-field}
		\end{equation}
		The field $K$ need not have finite degree over $\mathbb Q$.
	\end{lemma}
	\begin{proof}
		The roots $\lambda_i^n$ are distinct, so $g^n$ is
		regular semisimple. Put $L=K(\lambda_1,\ldots,\lambda_m)$
		and $L_n=K(\lambda_1^n,\ldots,\lambda_m^n)$.
		The extension $L/K$ is finite Galois, since it is a
		separable splitting field. Every
		$\sigma\in\operatorname{Gal}(L/L_n)$ permutes the
		$\lambda_i$ and fixes each $\lambda_i^n$. Injectivity
		of $\lambda_i\mapsto\lambda_i^n$ forces it to fix each
		$\lambda_i$. This group is therefore trivial, and the
		Galois correspondence gives $L_n=L$. The finiteness
		used here is that of $L/K$, not that of $K/\mathbb Q$.
	\end{proof}
	
	\begin{corollary}
		\label{cor:norm-realizing-powers}
		Under all the hypotheses of Theorem~\ref{thm:norm-numeric},
		one can choose a single $\gamma\in\Gamma$ such that,
		for every integer $n\ge1$, the power $\gamma^n$ is
		regular semisimple and loxodromic,
		\begin{equation}
			L_n=K(E(\gamma^n))=L=K(E(\gamma)),\qquad
			q_j\in N_{L/L^c}(L^\times)\quad(3\le j\le m-1).
			\label{eq:all-powers}
		\end{equation}
		The norm witnesses are the same for every $n$.
		Moreover, for each finite-index subgroup $\Delta\le\Gamma$,
		there is an $n\ge1$ with $\gamma^n\in\Delta$, without
		changing $L$, its actual complex conjugation, or these
		witnesses.
	\end{corollary}
	\begin{proof}
		Use the theorem's construction with the same $S$, $t$, and return
		sequence. Without residual factors, $\gamma=t$ already has every
		power regular by Proposition~\ref{prop:filtered-expansion}, excluding
		all torsion quotients of distinct roots.
		
		Otherwise, refine the targets at one residual factor using
		Lemma~\ref{lem:deep-targets}, retaining the norm conditions at all
		other factors. Corollary~\ref{cor:one-word} selects the connectors
		and then a common return exponent placing the final word in these
		smaller neighborhoods. In $B_c$, all roots lie in the deep unit group.
		In $B_u$, the central roots do so, whereas the extreme roots have
		absolute values greater and less than one. Every quotient involving
		an extreme root has absolute value different from one; quotients
		of distinct central roots lie in the torsion-free unit group.
		Thus no quotient of distinct roots is torsion.
		
		A field embedding preserves $z^n=1$, so this condition at one local
		realization excludes global torsion quotients. It is preserved by
		natural or dual transport and by extension from a source place.
		Lemma~\ref{lem:power-field} gives $L_n=L$ for every $n$.
		Positive powers remain loxodromic, and the same witnesses for
		$q_j=x_jc(x_j)$ remain in $L$ with the same actual complex conjugation.
		No intersection of infinitely many neighborhoods, or separate choice
		of exponent for each $n$, is needed for the powers.
		
		Finally, the action of $\Gamma$ on the finite coset set of $\Delta$
		has finite image. Some positive power of $\gamma$ acts trivially,
		hence belongs to $\Delta$. The exponent may depend on $\Delta$;
		no common power for every finite-index subgroup is asserted.
	\end{proof}
	
	\section{Algebraic fields and global unitary descent}
	\label{sec:global}
	
	Over algebraic fields, the number-field argument applies at a finite
	number-field stage, and its norm identities persist under scalar extension.
	One spectrum standardizes a normalized form; for an arbitrary initial form,
	a preliminary spectrum supplies a hyperbolic plane, giving the two-spectrum
	bound below.
	
	\subsection{Scalar extension and a finite dense core}
	
	\begin{lemma}
		\label{lem:global-fields}
		Let $K\subset\mathbb C$ be stable under actual complex
		conjugation $c$, with $[K:K^c]=2$, and let
		$g\in\mathrm{SU}(H)(K^c)$. Then $K(E(g))$ is stable
		under $c$. If $G/k$ is an affine algebraic group and
		$\Gamma\subset G(k)$ is Zariski-dense, then $\Gamma$ remains
		Zariski-dense after every field extension $k'/k$, viewed as a
		subset of $G_{k'}(k')$.
	\end{lemma}
	\begin{proof}
		The identity $g^{-1}=H^{-1}g^*H$ identifies the spectrum
		of $g^{-1}$ with that of $\bar g$. Hence
		$c(E(g))=E(g)^{-1}$, and the generated field is stable.
		For the second assertion, write a regular function on
		the extended group as a finite sum $\sum_i a_i f_i$,
		where the $a_i$ are linearly independent over $k$ and
		the $f_i$ are defined over $k$. If this function
		vanishes at every point of $\Gamma$, the values $f_i(\gamma)$ lie in $k$,
		so each $f_i$
		vanishes there and is zero by density. The argument
		takes place in the affine coordinate ring of the
		group and does not require the scalar extension
		to be algebraic.
	\end{proof}
	
	\begin{lemma}
		\label{lem:finite-core}
		Let $G$ be a connected algebraic group in
		characteristic zero with simple Lie algebra,
		and let $\Gamma\subset G(k)$ be a Zariski-dense
		subgroup. If $\Gamma$ contains an element of
		infinite order, it contains a finitely generated
		Zariski-dense subgroup.
	\end{lemma}
	\begin{proof}
		Let $A\in\Gamma$ have infinite order, and put
		$T=(\overline{\langle A\rangle}^{\rm Zar})^\circ$.
		Its dimension is positive, since a
		zero-dimensional algebraic group has only
		finitely many points over an algebraic
		closure. The sum
		\[
		\mathfrak a=\sum_{g\in\Gamma}\operatorname{Ad}(g)\operatorname{Lie}T
		\]
		is nonzero and stable under $\Gamma$. Its
		stabilizer is a closed algebraic subgroup,
		so density makes the subspace stable under $G$.
		It is therefore a nonzero ideal, hence equals
		$\operatorname{Lie}G$. A finite sum, with
		conjugating elements $g_1,\ldots,g_r$, already
		spans this space. The Zariski closure of
		$\Delta=\langle g_iAg_i^{-1}:1\leq i\leq r\rangle$
		contains every $g_iTg_i^{-1}$ and has the full
		Lie algebra of $G$. Since $G$ is connected,
		this closure is $G$.
	\end{proof}
	
	\subsection{Simultaneous norms over algebraic fields}
	
	\begin{proposition}
		\label{prop:algebraic-transfer}
		The conclusion of Theorem~\ref{thm:norm-numeric},
		with one element for all coefficients, together
		with every conclusion of
		Corollary~\ref{cor:norm-realizing-powers}, extends
		to every algebraic field
		$K\subset\overline{\mathbb Q}\subset\mathbb C$,
		possibly of infinite degree, that is stable
		under $c$ and satisfies $[K:K^c]=2$.
		The normalized model, the distinguished
		signature, and the prescribed group are
		retained.
	\end{proposition}
	\begin{proof}
		Put $k=K^c$. By Lemma~\ref{lem:global-fields}, density persists
		in the distinguished real realization. Real proximality
		\cite[\S\S3.1--3.2]{Benoist1997} gives a loxodromic element,
		hence one of infinite order. Since $\mathrm{SU}(H)$ is absolutely
		almost simple, Lemma~\ref{lem:finite-core} gives a finitely generated
		Zariski-dense subgroup $\Delta\subset\Gamma$.
		
		The entries of its generators, the entries
		of $H$, the coefficients $q_j$, a nonreal
		element of $K$, and all their complex
		conjugates under $c$ generate a number field
		$K_*\subset K$ stable under $c$.
		Put $k_*=K_*^c$. Then $[K_*:k_*]=2$ and
		$K=kK_*$: the compositum on the right is a
		subfield of $K$ not fixed by $c$ and already
		has degree two over $k$. The closure of
		$\Delta$ in the unitary model over $k_*$
		cannot be proper, since a proper closure
		would remain proper after the faithful
		extension of scalars to $k$, contradicting
		density.
		
		Apply Theorem~\ref{thm:norm-numeric} and
		Corollary~\ref{cor:norm-realizing-powers} to
		$\Delta$ over $K_*/k_*$. They give one
		$\gamma\in\Delta$ and elements
		$x_j\in L_*=K_*(E(\gamma))$ such that
		\[
		x_jc(x_j)=q_j \qquad(3\leq j\leq m-1).
		\]
		These identities remain valid in
		$L=KL_*=K(E(\gamma))$. The field $L$ is
		stable under $c$, and the involution remains
		nontrivial. Hence $[L:L^c]=2$, and the same
		identities give the required norms.
		For every $n\ge1$, equality of the
		number-field splitting fields gives
		\[
		K(E(\gamma^n))=K\,K_*(E(\gamma^n))
		=K\,K_*(E(\gamma))=L.
		\]
		Loxodromy and regularity of all powers persist under scalar extension.
		The same $x_j$ witness the norms for every $n$, and the finite-index
		conclusion uses the same action on cosets. The Hasse norm theorem is
		used only over $L_*/L_*^c$.
	\end{proof}
	
	\begin{remark}
		The proof uses simplicity and an element of
		infinite order. A Zariski-dense subgroup of
		an arbitrary algebraic group need not contain
		a finitely generated dense subgroup, as the
		roots of unity in a torus show. The argument
		also does not extend the result to
		transcendental fields: constructing $K_*$
		as a number field uses algebraicity.
	\end{remark}
	
	\subsection{Standard models from one or two spectra}
	
	Simultaneous norm witnesses give one common diagonal change of basis
	for the whole group. We first apply this to the normalized form;
	for an arbitrary initial form, a preliminary spectrum supplies
	the hyperbolic plane needed for normalization.
	
	\begin{corollary}
		\label{cor:normalized-descent}
		In the number-field setting of
		Theorem~\ref{thm:norm-numeric}, or the
		algebraic setting of
		Proposition~\ref{prop:algebraic-transfer},
		choose $\gamma$ with the additional
		properties of
		Corollary~\ref{cor:norm-realizing-powers},
		transferred to the algebraic setting when
		necessary. The individual splitting field
		$L=K(E(\gamma))$ defines a common change
		of basis realizing the entire prescribed
		group $\Gamma$ in the model
		\[
		J_m=\mathbb H\perp I_{m-2}.
		\]
		The same change of basis is valid over
		$K(E(\gamma^n))$ for every $n\ge1$ and
		over $F_{\rm eig}(\Gamma)$.
	\end{corollary}
	\begin{proof}
		Let $x_j\in L^\times$ be the simultaneous
		norm witnesses $q_j=x_jc(x_j)$. The matrix
		\[
		P=\operatorname{diag}(I_2,1,x_3^{-1},\ldots,x_{m-1}^{-1})
		\in\mathrm{GL}_m(L)
		\]
		satisfies $P^*HP=J_m$. Thus
		$P^{-1}\Gamma P\subset\mathrm{SU}(J_m)(L^c)$,
		with all entries in $L$ and determinant one.
		Constancy of the splitting fields of powers and the inclusion
		$L\subset F_{\rm eig}(\Gamma)$ give the remaining assertions.
	\end{proof}
	
	\begin{corollary}
		\label{cor:two-spectrum-algebraic}
		Let $K\subset\overline{\Q}\subset\C$ be stable under actual complex
		conjugation $c$, put $k=K^c$, and assume $[K:k]=2$. Let $H$ be a
		nonsingular $m$-dimensional Hermitian form over $K/k$, with $m\ge4$,
		of signature $(m-1,1)$ at the distinguished embedding. For every
		$k$-Zariski-dense subgroup $\Gamma\subset\SU(H)(k)$, there exist
		loxodromic elements $A,\gamma\in\Gamma$ and a finite $c$-stable extension
		\[
		L=K(E(A),E(\gamma))\subset F_{\rm eig}(\Gamma),
		\qquad [L:K]\le(m!)^2,
		\]
		over which $H$ is similar to $J_m=\mathbb H\perp I_{m-2}$.
		Consequently, one change of basis realizes the whole group in
		$\SU(J_m)(L^c)$. The element $\gamma$ can be chosen so that every
		positive power is regular semisimple and
		\[
		K(E(\gamma^n))=K(E(\gamma)),\qquad
		K(E(A),E(\gamma^n))=L\quad(n\ge1).
		\]
		Neither finite generation of $\Gamma$ nor isotropy over $K$ is assumed.
	\end{corollary}
	
	\begin{proof}
		By Lemma~\ref{lem:global-fields}, density persists in the distinguished
		real realization. Lemma~\ref{lem:proximal-semigroup}, applied in real
		rank one, gives a loxodromic $A\in\Gamma$. Put $K_1=K(E(A))$ and
		$k_1=K_1^c$. The field $K_1$ is algebraic over $\Q$ and stable under $c$.
		Since it contains $K$, conjugation is nontrivial on $K_1$ and
		$[K_1:k_1]=2$.
		
		The expanding and contracting roots of $A$ are simple. Their eigenlines
		have vectors over $K_1$, are isotropic, and pair nontrivially: otherwise
		they would span a totally isotropic plane in a Hermitian space of
		complex Witt index one. They therefore give a hyperbolic plane over
		$K_1$. Its orthogonal complement is positive definite at the
		distinguished embedding and can be diagonalized over $K_1/k_1$ as
		$\langle a_1,\ldots,a_{m-2}\rangle$, with $a_i\in k_1^\times$.
		Scaling the whole form by $a_1^{-1}$, and rescaling one vector in the
		hyperbolic pair, yields $P\in\GL_m(K_1)$ such that
		\[
		a_1^{-1}P^*HP
		=\mathbb H\perp\langle1,q_3,\ldots,q_{m-1}\rangle,
		\qquad q_j=a_{j-1}/a_1.
		\]
		Thus $P^{-1}\Gamma P$ is a prescribed $k_1$-Zariski-dense subgroup of
		this normalized model. Density follows by scalar extension in
		Lemma~\ref{lem:global-fields}, followed by the isomorphism induced by
		$P$.
		
		Apply Proposition~\ref{prop:algebraic-transfer} and
		Corollary~\ref{cor:normalized-descent} over $K_1/k_1$. They give an
		element $P^{-1}\gamma P$, with $\gamma\in\Gamma$, and simultaneous
		norm witnesses over $K_1(E(\gamma))$. The resulting diagonal change
		of basis standardizes the normalized form. Composing it with $P$
		proves the asserted similarity over
		$L=K(E(A),E(\gamma))$. Conjugacy does not change either spectrum.
		Both spectra belong to matrices over $K$ of degree $m$, so each
		splitting field over $K$ has degree at most $m!$, giving the displayed
		bound.
		
		Every positive power of the selected $\gamma$ is regular. Equivalently,
		no quotient of two distinct roots is a root of unity. Applying
		Lemma~\ref{lem:power-field} directly over the original field $K$ gives
		the stated equality of individual splitting fields; adjoining $E(A)$
		gives the second equality. No regularity assumption on the central
		roots of the preliminary element $A$ is needed.
	\end{proof}

	\begin{proposition}
		\label{prop:global-transfer}
		Let $K/\Q$ be imaginary quadratic, let $H_K$ be a Hermitian form over
		$K/\Q$ of signature $(3,1)$,
		and let $\Gamma\subset\SU(H_K)(\Q)$ be Zariski-dense. There are
		$A,\gamma\in\Gamma$ and a finite $c$-stable field
		\[
		L=K(E(A),E(\gamma))\subset F_{\rm eig}(\Gamma),\qquad [L:K]\le576,
		\]
		over which $H_K$ is similar to $J=\mathbb H\perp\langle1,1\rangle$.
		One change of basis realizes all of $\Gamma$ in $\SU(J)(L^c)$.
	\end{proposition}
	\begin{proof}
		Apply Corollary~\ref{cor:two-spectrum-algebraic} with $m=4$ and
		$k=\Q$. The degree bound is $(4!)^2=576$. In particular, its conclusion
		holds for the finitely generated groups in
		\cite[Conjecture~1.5]{CanoCastroGusevskii}.
	\end{proof}
	
	\begin{remark}
		The initial spectrum supplies a hyperbolic plane over $K(E(A))$;
		the second makes the normalized coefficients norms. Isotropy over
		the initial field is not assumed. The argument gives two spectra
		over the original $K$, with no assertion that one spectrum always
		suffices for an arbitrary initial form. Finite generation is unnecessary
		for this normalization and for the resulting descent.
	\end{remark}

	\section{A transcendental obstruction in intrinsic and total spectral fields}
	\label{sec:transcendental}
	
	We retain the conventions of Section~\ref{sec:introduction}, with $c$
	denoting actual complex conjugation in the fixed embedding into $\C$.
	A loxodromic element of $\SU(3,1)$ has one simple expanding root,
	one simple contracting root, and two central roots of modulus one;
	the element is regular precisely when the central roots are distinct. The construction
	below concerns intrinsic fields of positive transcendence degree.
	A $c$-invariant valuation gives every norm even value and the residual
	Hermitian coefficient odd value. The argument does not use the
	number-field results or the Hasse norm theorem.
	
	\begin{theorem}
		\label{thm:transcendental-obstruction}
		There exists a finitely generated subgroup $\Gamma\subset\SU(3,1)$
		that is free of rank four, discrete, torsion-free, convex cocompact,
		and Zariski-dense, with the following properties.
		\begin{enumerate}[label=\textup{(\roman*)},leftmargin=*]
			\item Every nontrivial element of $\Gamma$ is loxodromic and regular
			semisimple, and all its nontrivial powers are regular.
			\item For every $A\in\Gamma\setminus\{1\}$, the intrinsic field
			\[
			K_A=\mathbb Q(\tr\Gamma,E(A))\subset\mathbb C,
			\qquad k_A=K_A^c
			\]
			is finitely generated over $\mathbb Q$, has positive transcendence
			degree, and satisfies $[K_A:k_A]=2$. The involution $c$ is actual
			complex conjugation in the fixed embedding.
			\item Matrix extraction over $K_A$, with $A$ diagonal, yields a
			conjugate representation preserving
			\[
			H_A=\mathbb H\perp\langle1,q_A\rangle,
			\qquad q_A\in k_A^\times,\quad q_A>0.
			\]
			Fix any such normalized extraction before choosing further elements
			whose spectra are to be adjoined. For every finite family
			$\gamma_1,\ldots,\gamma_r\in\Gamma\setminus\{1\}$, including $r=0$,
			put
			$L_{A,\boldsymbol\gamma}=K_A(E(\gamma_1),\ldots,E(\gamma_r))$,
			with $L_{A,\boldsymbol\gamma}=K_A$ when $r=0$. Then
			\begin{equation}\label{eq:finite-spectral-nonnorm}
				q_A\notin N_{L_{A,\boldsymbol\gamma}/L_{A,\boldsymbol\gamma}^c}
				(L_{A,\boldsymbol\gamma}^{\times}).
			\end{equation}
			\item The obstruction persists in the intrinsic total spectral field:
			\begin{equation}\label{eq:total-spectral-nonnorm}
				q_A\notin N_{L_{A,\infty}/L_{A,\infty}^c}
				(L_{A,\infty}^{\times}),
				\qquad L_{A,\infty}=K_A(E(\gamma):\gamma\in\Gamma).
			\end{equation}
			This total field is independent of $A$ and equals
			$L_{\rm spec}(\Gamma)=\Q(E(g):g\in\Gamma)$.
			Moreover, $q_A$ itself is transcendental over $\mathbb Q$.
		\end{enumerate}
		Thus the universal spectral norm assertion for intrinsic fields of
		positive transcendence degree fails already in dimension four. The
		group has no nontrivial unipotent elements. The obstruction is present
		before any word is chosen and persists for every choice of $A$.
	\end{theorem}
	
	\paragraph{Order of choices.}
	Fix the four generators, the common interval, and a transcendental
	parameter $T$; these determine a single group. Choose any $A\ne1$
	and any normalized extraction over $K_A$, and keep its coefficient
	$q_A$ fixed. The obstruction must then hold for every finite family
	of nontrivial words, including the empty family. The auxiliary valuation
	may depend on that family; the group and the extraction do not depend
	on a putative norm witness.
	
	\subsection{A uniform family of representations}
	
	\begin{lemma}
		\label{lem:uniform-seed}
		Let
		\[
		J_3=\begin{pmatrix}0&1&0\\1&0&0\\0&0&1\end{pmatrix},
		\qquad G_0=\SU(J_3)(\mathbb R),
		\]
		where the notation for real points refers to the unitary group over
		$\mathbb C/\mathbb R$. There exist $B_1,\ldots,B_4\in G_0$ such that:
		\begin{enumerate}[label=\textup{(\roman*)}]
			\item the representation $\rho_B:F_4\to G_0$ sending the free
			generators to $B_1,\ldots,B_4$ is faithful and convex cocompact;
			\item each of the pairs $(B_1,B_2)$ and $(B_3,B_4)$ generates a
			Zariski-dense subgroup of $G_0$;
			\item for every nontrivial word $w\in F_4$, the matrix
			$B_w=\rho_B(w)$ satisfies $\det(B_w-I_3)\ne0$.
		\end{enumerate}
		In particular, for every nontrivial word $w$, the polynomial
		$(X-1)\chi_{B_w}(X)$ has four distinct roots: one of modulus greater
		than one, one of modulus less than one, and two of modulus one.
		These assertions also hold for all nontrivial powers.
	\end{lemma}
	
	\begin{proof}
		Choose four hyperbolic elements of $\SO^+(2,1)$ with eight distinct
		fixed points and take powers satisfying strict ping-pong inclusions.
		This gives a faithful convex cocompact representation of $F_4$;
		the totally geodesic embedding into the complex hyperbolic plane
		preserves its convex hull and compact quotient. See also
		\cite[\S6.1\textup{(ii)}]{GuichardWienhard2012}.
		
		By \cite[Theorems~5.13 and~5.15]{GuichardWienhard2012}, the
		Anosov property is open in $\operatorname{Hom}(F_4,G_0)=G_0^4$
		and, in real rank one, is equivalent to finite kernel and convex
		cocompact image. The finite kernel is trivial since $F_4$ is torsion-free.
		Thus there is a nonempty open set $\mathcal S$ of faithful convex
		cocompact representations.
		
		For each $w\ne1$, the function
		\[
		f_w(B_1,\ldots,B_4)=\det\bigl(w(B_1,\ldots,B_4)-I_3\bigr)
		\]
		is not identically zero. Indeed, the complexification of the real
		algebraic group $G_0$ is $\SL_{3,\mathbb C}$. Borel's word-map
		theorem \cite[\S1, Theorem~1, p.~154]{Borel1983} states that the
		map defined by a nontrivial word is dominant on this group. The
		function $g\mapsto\det(g-I_3)$ does not vanish identically on
		$\SL_3$, as witnessed by $\operatorname{diag}(2,3,1/6)$.
		Its composition with the word map is therefore nonzero. Since the
		real form is Zariski-dense in its complexification, $f_w$ cannot
		vanish on all of $G_0^4$.
		
		Consequently, $|f_w|^2$ is a nonzero real analytic function on the
		connected real analytic manifold $G_0^4$. Its zero set has Haar
		measure zero. There are only countably many nontrivial words in
		$F_4$, so the union of these exceptional sets still has measure zero.
		
		By \cite[Proposition~16]{Winkelmann2002}, the pairs that do not
		generate a Zariski-dense subgroup form a measure-zero subset of
		$G_0^2$. Their inverse images under the projections to factors
		$(1,2)$ and $(3,4)$ have measure zero in $G_0^4$. This follows from
		Fubini's theorem, since Haar measure on a Lie group is
		$\sigma$-finite. A nonempty relatively compact open subset of
		$\mathcal S$ has positive finite measure and cannot be covered by
		the countable union of all the exceptional sets above. Choosing a
		quadruple in its complement gives \textup{(i)}--\textup{(iii)}
		simultaneously.
		
		Every nontrivial $B_w$ is loxodromic: a faithful convex cocompact
		representation of a free group has neither torsion nor parabolic
		elements, and discreteness excludes elliptic elements of infinite
		order. In dimension three, its three roots have moduli $>1$, $<1$,
		and $1$, respectively, and hence are distinct. Condition
		\textup{(iii)} separates the additional root $1$ from them.
		Finally, $w^r\ne1$ whenever $w\ne1$ and $r\ne0$, so the same
		argument applies to all nontrivial powers with the same generators.
	\end{proof}
	
	\begin{proposition}
		\label{lem:unitary-deformation}
		Let $B_1,\ldots,B_4$ be as in the preceding lemma, and let
		$F\subset\mathbb R$ be the field generated over $\mathbb Q$ by the
		real and imaginary parts of their entries. For an indeterminate $t$,
		put
		\[
		H_t=\operatorname{diag}(J_3,t),\qquad
		M(B)=\operatorname{diag}(B,1),\qquad
		R_t=I_2\oplus\frac1{1+t}
		\begin{pmatrix}1-t&-2t\\2&1-t\end{pmatrix},
		\]
		and define
		\[
		g_1(t)=M(B_1),\quad g_2(t)=M(B_2),\quad
		g_j(t)=R_tM(B_j)R_t^{-1}\quad(j=3,4).
		\]
		There exists $\delta>0$ such that, for every $0<t<\delta$, these
		four generators define a faithful, convex cocompact, Zariski-dense
		representation of $F_4$ in $\SU(H_t)$. In particular, every nontrivial
		word is loxodromic for every parameter in this common interval.
		
		If $T\in(0,\min\{1,\delta\})$ is transcendental over $F$, every
		nontrivial word in the specialized generators, including every
		nontrivial power, is also regular semisimple.
	\end{proposition}
	
	\begin{proof}
		Set $a=(1-t)/(1+t)$ and $b=2/(1+t)$. The identity $a^2+tb^2=1$
		gives $R_t^*H_tR_t=H_t$ and $\det R_t=1$. A direct calculation
		with matrix units also gives, for $j=1,2$,
		\begin{equation}\label{eq:deformation-lie}
			\operatorname{Ad}(R_t)(E_{j3})=aE_{j3}+tbE_{j4},\qquad
			\operatorname{Ad}(R_t)(E_{3j})=aE_{3j}+bE_{4j}.
		\end{equation}
		At $t=0$,
		\begin{equation}\label{eq:triangular-reduction}
			R_0=\begin{pmatrix}I_3&0\\2e_3^t&1\end{pmatrix},\qquad
			R_0M(B)R_0^{-1}=\begin{pmatrix}B&0\\2e_3^t(B-I_3)&1\end{pmatrix}.
		\end{equation}
		Matrices of this block triangular form constitute a group, and
		projection to the upper left $3\times3$ block is a homomorphism.
		Thus every word, including words with inverses, specializes at zero
		to a matrix with upper left block $B_w$ and final diagonal entry
		$1$. In particular,
		$\chi_{g_w(0)}(X)=(X-1)\chi_{B_w}(X)$, with no restriction on
		word length.
		
		For $t>0$, let $D_t=\operatorname{diag}(1,1,1,\sqrt t)$, so that
		$D_t^*H_1D_t=H_t$. Conjugation by $D_t$ puts all the groups in
		the fixed group $G_1=\SU(H_1)(\mathbb R)$. Since
		\[
		D_tR_tD_t^{-1}
		=I_2\oplus\frac1{1+t}
		\begin{pmatrix}1-t&-2\sqrt t\\2\sqrt t&1-t\end{pmatrix}
		\longrightarrow I_4,
		\]
		the resulting representations, viewed as points of
		$\operatorname{Hom}(F_4,G_1)=G_1^4$, converge to
		$w\mapsto M(B_w)$ as $t\downarrow0$.
		
		The limiting representation is convex cocompact in $G_1$. Indeed,
		the copy of $\mathbb H^2_{\mathbb C}$ defined by the first three
		coordinates is a closed, totally geodesic, convex subspace of
		$\mathbb H^3_{\mathbb C}$. The convex hull of the limit set lies
		in this subspace and is the same whether computed there or in the
		ambient space; its quotient remains compact. No irreducibility or
		Zariski-density assumption on the limiting representation in $G_1$
		is needed.
		
		Apply \cite[Theorems~5.13 and~5.15]{GuichardWienhard2012} to the
		semisimple real rank-one group $G_1$, which has finite center, and
		to the finitely generated word-hyperbolic group $F_4$. After
		decreasing $\delta$ if necessary, a single neighborhood of the
		limiting representation contains the entire family with
		$0<t<\delta$. Every representation in this neighborhood has
		convex cocompact image and finite, hence trivial, kernel. This
		proves the assertion for all words simultaneously: the interval
		comes from an open set of representations, rather than from an
		intersection of intervals chosen separately for individual words.
		
		Fix any $t\in(0,\delta)$. The two dense pairs in the preceding
		lemma ensure that the Zariski closure contains both
		$M(\SU(J_3))$ and $R_tM(\SU(J_3))R_t^{-1}$. Over an algebraic
		closure, its Lie algebra contains $\mathfrak{sl}_3\oplus0$ and
		its conjugate by $R_t$. In \eqref{eq:deformation-lie}, both
		$b=2/(1+t)$ and $tb$ are nonzero. Subtracting the terms already
		in $\mathfrak{sl}_3$ gives $E_{j4}$ and $E_{4j}$ for $j=1,2$.
		Their brackets with $\mathfrak{sl}_3$ give $E_{34}$ and $E_{43}$.
		In addition, $[E_{j4},E_{4j}]=E_{jj}-E_{44}$ supplies the remaining
		diagonal direction. Thus the Lie algebra is
		$\mathfrak{sl}_4$. Connectedness of the ambient group implies
		Zariski density for every parameter in the common interval;
		this step does not require transcendence.
		
		For each $w\ne1$, \eqref{eq:triangular-reduction} gives
		\[
		\chi_{g_w(t)}(X)\bmod t=(X-1)\chi_{B_w}(X).
		\]
		The right-hand side is separable by the lemma. Hence the
		discriminant of the characteristic polynomial on the left is a
		rational function in $F(i)(t)$ that is regular and nonzero at
		$t=0$. It is also nonzero at $t=T$, because $T$ is transcendental
		over $F(i)$. Thus every nontrivial word, including all its nontrivial powers,
		is regular at the transcendental specialization $T$.
	\end{proof}
	
	\subsection{Actual complex conjugation and finite composita}
	
	To obstruct a norm after adjoining several spectra, we need a
	valuation on their common field that is invariant under the
	conjugation inherited from $\C$. The next lemma constructs
	this valuation by continuing all the root sets along one common path.
	
	\begin{lemma}
		\label{lem:conjugation-compatible-lifting}
		Let $E\subset\mathbb C$ be a field stable under complex conjugation,
		with $i\in E$, and let $E'/E$ be a finite extension inside
		$\mathbb C$ that is also stable under complex conjugation. Let
		$P_1(t,X),\ldots,P_r(t,X)\in E(t)[X]$ be monic polynomials of
		degree four with constant term in $X$ equal to one. Suppose that:
		\begin{enumerate}[label=\textup{(\roman*)},leftmargin=*]
			\item their coefficients are regular at $t=0$ and on an interval
			$I=(0,\delta)$, and they satisfy
			\[
			\overline{P_j}(t,X)=X^4P_j(t,X^{-1});
			\]
			\item each $P_j(0,X)$ has four distinct roots: one of modulus
			greater than one, one of modulus less than one, and two of
			modulus one;
			\item for every $t\in I$, each $P_j(t,X)$ has exactly one simple
			expanding root, one simple contracting root, and two roots of
			modulus one, counted with multiplicity.
		\end{enumerate}
		Fix $T\in I$ transcendental over $E'$, and put
		\[
		M=E'(T),\qquad
		L=M\bigl(\text{all roots of }P_j(T,X),\ 1\le j\le r\bigr)
		\subset\mathbb C.
		\]
		Then $L$ is stable under actual complex conjugation $c$ and admits
		a discrete valuation $w:L^\times\to\mathbb Z$ such that
		\[
		w(T)=1,\qquad w(c(z))=w(z)\quad(z\in L^\times).
		\]
		The valuation is trivial on the constant field: $w(a)=0$ for every
		$a\in(E')^\times$.
		In particular, $T\notin N_{L/L^c}(L^\times)$.
		The product $\prod_jP_j(0,X)$ need not be separable.
	\end{lemma}
	
	\begin{proof}
		\emph{Convergent series and the formal involution.}
		Let $E_0\subset\mathbb C$ be the compositum of the splitting fields
		of the polynomials $P_j(0,X)$ over $E'$. This is a finite extension
		of $E'$ and is stable under actual complex conjugation, which sends
		each root set to its set of inverses. Label the roots of $P_j(0,X)$
		as $\beta_{j,+},\beta_{j,-},\beta_{j,1},\beta_{j,2}$ according to
		their moduli. Hensel's lemma in the complete discrete valuation ring
		$E_0[[t]]$ \cite[Proposition~7.31]{MilneANT2020} gives, for each
		root $\beta$, a unique series
		\[
		\Lambda_{j,\beta}(t)\in E_0[[t]],\qquad
		P_j(t,\Lambda_{j,\beta}(t))=0,\qquad
		\Lambda_{j,\beta}(0)=\beta.
		\]
		These series converge as complex functions near zero, by the
		holomorphic implicit function theorem at a simple root. We apply
		Hensel's lemma to \emph{each} polynomial separately; shared roots
		of two different polynomials at $t=0$ do not affect the argument.
		
		Coefficientwise conjugation defines an involution $c_0$ of
		$E_0((t))$ fixing $t$. The reciprocal identity and uniqueness of
		the lifts give
		\[
		c_0(\Lambda_{j,+})=\Lambda_{j,-}^{-1},\qquad
		c_0(\Lambda_{j,-})=\Lambda_{j,+}^{-1},\qquad
		c_0(\Lambda_{j,a})=\Lambda_{j,a}^{-1}\quad(a=1,2).
		\]
		Indeed, $c_0(\Lambda_{j,\beta})^{-1}$ is the root of $P_j$ whose
		constant term is $\bar\beta^{-1}$. For a root of modulus one this
		constant term is $\beta$, whereas the expanding and contracting
		roots are exchanged. Hence the field of algebraic germs
		\[
		\mathcal L=E'(t)(\Lambda_{j,+},\Lambda_{j,-},
		\Lambda_{j,1},\Lambda_{j,2}:1\le j\le r)
		\subset E_0((t))
		\]
		is stable under $c_0$. The restriction of $\operatorname{ord}_t$
		to $\mathcal L$ is $c_0$-invariant and takes the value one on $t$.
		It remains to identify $c_0$ with actual complex conjugation on $L$.
		
		\emph{A common path and control of collisions.}
		Each discriminant $\Delta_j(t)$ is a nonzero rational function,
		since $\Delta_j(0)\ne0$. It has only finitely many zeros, and
		$\Delta_j(T)\ne0$ because $T$ is transcendental over $E'$.
		Choose $0<\epsilon<T$ sufficiently small that all the germs above
		are regular on $[0,\epsilon]$. The real zeros of the individual
		discriminants in $[\epsilon,T]$ form a finite set. By assumption,
		none of the coefficients has a pole on this segment.
		
		At any real zero $a$ of one of these discriminants, the expanding
		and contracting roots of \emph{every} polynomial remain simple.
		Both extend holomorphically to a complex disc around $a$. After
		shrinking the disc, their values remain separated from each other
		and from the other two roots. This follows from continuity of
		roots, using disjoint contours around their values at $a$. Any
		collision at $a$ therefore involves only the two central roots of
		the relevant polynomial. Since the family is finite, a common disc
		can be chosen for all its members.
		
		Choose these discs to be disjoint and sufficiently small to contain
		no other discriminant zeros or coefficient poles. Replace the
		passage through each real center by a small semicircle inside its
		disc. This gives a single path from $\epsilon$ to $T$ avoiding the
		union of the zeros of the \emph{individual} discriminants. On each
		real segment, continuation preserves the unique expanding root
		and the unique contracting root. Along a detour, their holomorphic
		extensions to the whole disc prevent either from being exchanged
		with a central root. The central roots may be exchanged with each
		other, which will cause no difficulty. We do not use the
		discriminant of $\prod_jP_j$, which may vanish even at zero.
		
		\emph{Joint continuation preserves field relations.}
		Continue $\mathcal L$ as a field of meromorphic germs along this
		same path, using joint continuation for all its generators. This
		preserves all algebraic relations among the generators; the root
		bijections at the endpoint are not chosen independently.
		Evaluation at $T$ defines an isomorphism
		\[
		\Psi:\mathcal L\longrightarrow L,
		\qquad \Psi(t)=T,\quad \Psi|_{E'}=\operatorname{id}.
		\]
		To justify this, we must exclude new relations and vanishing
		denominators under evaluation. In the ring generated by the roots
		over $E'(t)$, the root generators are analytic at the endpoint of
		the path. Every rational function in $E'(t)$ is also regular at
		$T$, since $T$ is transcendental over $E'$. If $b$ is a nonzero
		element of this ring, all its conjugates over $E'(t)$ are likewise
		polynomial expressions in the roots and are regular at the
		endpoint. Their product, counted with the multiplicities in the
		field norm, is a nonzero rational function in $E'(t)$, so its
		value at $T$ is nonzero.
		None of the conjugates, including $b$, can therefore vanish at
		$T$. Every quotient can consequently be evaluated, and evaluation
		is injective. Its image is precisely the field generated by all
		the root sets in the statement, proving surjectivity. This also
		shows that the finite enlargement of the constant field from $E$
		to $E'$ requires no number-field hypothesis.
		
		\emph{Comparison with actual complex conjugation.}
		Write $\lambda_{j,*}=\Psi(\Lambda_{j,*})$. By the control of the
		path, $\lambda_{j,+}$ and $\lambda_{j,-}$ are the expanding and
		contracting roots at $T$, and both $\lambda_{j,a}$ have modulus
		one. The reciprocal identity and uniqueness of the expanding and
		contracting roots imply
		\[
		c(\lambda_{j,+})=\lambda_{j,-}^{-1},\qquad
		c(\lambda_{j,-})=\lambda_{j,+}^{-1},\qquad
		c(\lambda_{j,a})=\lambda_{j,a}^{-1}.
		\]
		These are exactly the images of the identities for $c_0$. They
		remain valid if the central branches were exchanged along a
		detour. The involutions also correspond on $E'(t)$ because $T$
		is real and conjugation on $E'$ is actual complex conjugation.
		Since these elements generate $\mathcal L$, we obtain
		\[
		\Psi\circ c_0=c\circ\Psi.
		\]
		Thus $w(z)=\operatorname{ord}_t(\Psi^{-1}z)$ is a $c$-invariant
		valuation on $L$. Its value group is contained in $\mathbb Z$ and
		contains $1=w(T)$, so it is exactly $\mathbb Z$. Finally,
		$T=zc(z)$ would imply $1=w(T)=2w(z)$, a contradiction.
	\end{proof}
	
	\begin{lemma}
		\label{prop:finite-level-obstruction}
		Let $K\subset\mathbb C$ be stable under complex conjugation, with
		$c|_K\ne1$, and let $L_\infty=\bigcup_\alpha L_\alpha$ be a
		directed union of $c$-stable subfields containing $K$. For every
		nonzero $a\in K^c$,
		\[
		a\in N_{L_\infty/L_\infty^c}(L_\infty^\times)
		\quad\Longleftrightarrow\quad
		\exists\alpha:\ a\in N_{L_\alpha/L_\alpha^c}
		(L_\alpha^\times).
		\]
	\end{lemma}
	
	\begin{proof}
		Each field has degree two over its fixed field because $c$ is
		already nontrivial on $K$, and its norm is $z\mapsto zc(z)$.
		If $a=zc(z)$ in $L_\infty$, then $z$ belongs to some $L_\alpha$,
		as does $c(z)$ by stability. The same identity realizes $a$ as a
		norm at that stage. The converse is immediate. No compatibility
		between valuations chosen at different stages is required.
	\end{proof}
	
	\subsection{Extraction over the trace field with one spectrum adjoined}
	
	The deformation is expressed over a field containing the matrix
	entries. We now realize the same representation over the intrinsic
	field $K_A$ and compare the discriminant classes of the ambient
	and extracted Hermitian forms.
	
	\begin{lemma}
		\label{lem:trace-algebra}
		Let $M$ be a field of characteristic zero, and let
		$\Lambda\subset\SL_4(M)$ be a group whose $M$-linear span is
		$M_4(M)$. Put
		\[
		t=\mathbb Q(\tr\lambda:\lambda\in\Lambda).
		\]
		If $A\in\Lambda$ has four distinct roots in $M$ and
		$K=t(E(A))$, there exist $P\in\GL_4(M)$ and a representation
		\begin{equation}\label{eq:trace-field-realization}
			\rho_K(\lambda)=P^{-1}\lambda P\in\SL_4(K)
			\qquad(\lambda\in\Lambda).
		\end{equation}
		If $\Lambda$ is finitely generated, then both $t$ and $K$ are
		finitely generated over $\mathbb Q$.
	\end{lemma}
	
	\begin{proof}
		Choose an $M$-basis $B_1=I,B_2,\ldots,B_{16}$ of $M_4(M)$
		consisting of elements of $\Lambda$. The matrix
		\[
		G=(\tr(B_iB_j))_{i,j}
		\]
		is invertible because the pairing $(X,Y)\mapsto\tr(XY)$ on
		$M_4(M)$ is nondegenerate, as is seen using the matrix units.
		Its entries lie in $t$. The coordinates of each $\lambda\in\Lambda$
		in this basis are obtained from $G$ and the values
		$\tr(\lambda B_j)\in t$, so they too lie in $t$. Consequently,
		\[
		\mathcal A=\operatorname{span}_t\Lambda
		=\bigoplus_{i=1}^{16}tB_i
		\]
		is a $16$-dimensional algebra with
		$\mathcal A\otimes_tM=M_4(M)$. It is simple: a nonzero two-sided
		ideal would extend to a nonzero ideal of $M_4(M)$ and hence would
		already have dimension $16$. Its center is $t$, since after
		extension to $M$ every central element is a scalar matrix $zI$,
		and $4z=\tr(zI)\in t$.
		
		Write the roots of $A$ as $\lambda_1,\ldots,\lambda_4$. The
		algebra $\mathcal A_K=\mathcal A\otimes_tK$ contains the idempotent
		\[
		e=\prod_{j=2}^{4}
		\frac{A-\lambda_jI}{\lambda_1-\lambda_j}.
		\]
		After extension to $M$, this is the rank-one spectral projector
		corresponding to $\lambda_1$. The left ideal $V=\mathcal A_K e$
		has dimension $4$ over $K$, because
		$V\otimes_KM=M_4(M)e$ has that dimension over $M$. Left
		multiplication gives an injective map
		\[
		\mathcal A_K\longrightarrow\operatorname{End}_K(V).
		\]
		Injectivity follows from the simplicity of $\mathcal A_K$, proved
		by the same scalar-extension argument. Both algebras have dimension
		$16$, so this map is an isomorphism. This is the realization of a
		simple algebra on a minimal left ideal underlying the usual
		structure theorem; see \cite[Chapter~IV, Theorem~1.17]{MilneCFT2020}.
		
		To obtain the conjugating matrix over the stated field, choose
		$0\ne v\in eM^4$. The map $M_4(M)e\to M^4$ defined by
		$Xe\mapsto Xv$ is an isomorphism of modules. The images of a
		$K$-basis of $V$ form the columns of a matrix $P\in\GL_4(M)$
		and satisfy $\lambda P=P\rho_K(\lambda)$. This proves
		\eqref{eq:trace-field-realization}, including preservation of determinants
		and characteristic polynomials.
		
		Now suppose that $s_1,\ldots,s_r$ generate $\Lambda$. Let $t_0$
		be the field generated by the finite lists
		\[
		\tr(B_iB_j),\qquad \tr(B_iB_jB_k),\qquad
		\tr(s_\ell^{\pm1}B_j).
		\]
		Each entry is the trace of an element of $\Lambda$, so $t_0\subset t$.
		Inverting $G$ shows that the structure constants of the products
		$B_iB_j$, as well as the coordinates of every $s_\ell^{\pm1}$,
		belong to $t_0$. Hence $\bigoplus t_0B_i$ is an algebra containing
		$\Lambda$. Since $\tr(B_i)=\tr(B_iB_1)\in t_0$, every trace of
		an element of $\Lambda$ belongs to $t_0$, and therefore $t=t_0$.
		Newton's identities give $\chi_A\in t[X]$. Its splitting field
		$K$ has degree at most $4!$ over $t$, proving both finite-generation
		assertions.
	\end{proof}
	
	\begin{proposition}
		\label{prop:intrinsic-discriminant}
		Let $M\subset\mathbb C$ be stable under actual complex conjugation
		$c$, with $c|_M\ne1$, and put $M_0=M^c$. Let
		$\tau\in M_0^\times$ be positive in the distinguished embedding,
		and consider
		\[
		H_\tau=\mathbb H\perp\langle1,\tau\rangle.
		\]
		Let $\Lambda$ be a Zariski-dense subgroup of
		$\SU(H_\tau)(M_0)$, and let $A\in\Lambda$ be regular loxodromic,
		with $E(A)\subset M$. For the intrinsic field
		\[
		K_A=\mathbb Q(\tr\Lambda,E(A))\subset M,\qquad k_A=K_A^c,
		\]
		the following assertions hold.
		\begin{enumerate}[label=\textup{(\roman*)}]
			\item $c(K_A)=K_A$ and $[K_A:k_A]=2$.
			\item The same representation has a model over $K_A$, with $A$
			diagonal, whose image is $k_A$-Zariski-dense and preserves the
			normalized form
			\[
			H_A=\mathbb H\perp\langle1,q_A\rangle,\qquad
			q_A\in k_A^\times,\quad q_A>0
			\]
			in the distinguished embedding.
			\item There exist an intertwiner $P\in\GL_4(M)$ and a scalar
			$a\in M_0^\times$ such that $P^*H_\tau P=aH_A$.
			Consequently, the coefficient obtained by this extraction
			satisfies, with $u=\det P/a^2\in M^\times$,
			\begin{equation}\label{eq:discriminant-transfer}
				q_A=\tau\,u\,c(u)
				=\tau\,N_{M/M_0}\!\left(\frac{\det P}{a^2}\right).
			\end{equation}
		\end{enumerate}
		The coefficient $q_A$ is fixed by $A$ and the extraction before
		any further word is considered. Identity \eqref{eq:discriminant-transfer}
		holds for every extraction of this representation to the indicated
		normalized form. No assumption is made that $M$ is a function
		field or is algebraic over $\mathbb Q$, or that the spectrum of
		$A$ is constant in a family.
	\end{proposition}
	
	\begin{proof}
		By Lemma~\ref{lem:global-fields}, Zariski density persists under
		extension of scalars.
		After the unitary group is split over $M$, the Zariski closure
		of $\Lambda$ is $\SL_4$. Its linear span is $M_4(M)$:
		transvections give the off-diagonal matrix units, and their
		products give the diagonal ones. Lemma~\ref{lem:trace-algebra} therefore
		applies.
		
		Unitarity gives $c(\tr\lambda)=\tr(\lambda^{-1})$, so the trace
		field is stable under $c$. The field generated by the roots of
		$A$ is also stable, since conjugation sends each root to the
		inverse of another. Moreover, $c$ cannot act trivially on the
		trace field. Otherwise $\tr\lambda=\tr(\lambda^{-1})$ for every
		$\lambda\in\Lambda$, and Zariski density would extend this
		identity to $\SL_4$. But for $x=\diag(2,3,5,1/30)$,
		$\tr x-\tr(x^{-1})=-21\ne0$. Thus $c|_{K_A}$ is a nontrivial
		involution, and its fixed field has index two.
		
		The preceding lemma gives $P_0\in\GL_4(M)$ such that
		$\rho(\Lambda)=P_0^{-1}\Lambda P_0\subset\SL_4(K_A)$.
		The form $H'=P_0^*H_\tau P_0$ is Hermitian and invariant.
		The equations $\rho(\lambda)^*X\rho(\lambda)=X$ are linear
		over $K_A$. A finite subsystem has the same solution space,
		even when $\Lambda$ is infinite, since there are only sixteen
		unknowns. After extension to $M$, this solution space is a line:
		for any invariant matrix $X$, the matrix $(H')^{-1}X$ commutes
		with $\rho(\Lambda)$, whose linear span is the full matrix
		algebra, and hence is scalar. Descent of this linear system
		therefore gives a nonzero invariant matrix $J\in M_4(K_A)$.
		It is necessarily invertible, and the invariant matrices over
		$K_A$ form the line $K_AJ$.
		
		The adjoint $J^*=c(J)^t$ is also invariant, so write $J^*=zJ$
		with $z\in K_A^\times$. Applying the adjoint twice gives
		$zc(z)=1$. The quadratic case of Hilbert's theorem~90
		\cite[Chapter~II, Corollary~1.23]{MilneCFT2020} can here be
		implemented explicitly. If $z\ne-1$, set $a=1+z$; then
		$a/c(a)=z$ and $aJ$ is Hermitian. If $z=-1$, choose
		$0\ne a\in K_A$ with $c(a)=-a$; again $aJ$ is Hermitian.
		We have thus obtained an invariant Hermitian form $H_K$ over
		$K_A$. The Hermitian forms $H'$ and $H_K$ differ by a scalar
		$b\in M_0^\times$. Replacing $H_K$ by $-H_K$ if necessary,
		we may assume that its signature in the distinguished embedding
		is $(3,1)$.
		
		All four roots of $\rho(A)$ lie in $K_A$, so its four eigenlines
		have bases over $K_A$. The expanding and contracting eigenlines
		are isotropic and pair nontrivially with each other. The two
		central eigenlines are orthogonal to them and to each other.
		These assertions follow directly from
		\[
		(c(\lambda_i)\lambda_j-1)\,H_K(v_i,v_j)=0
		\]
		and simplicity of the spectrum. The coefficients on the central
		lines are $r,s\in k_A^\times$, both positive in the distinguished
		embedding. Rescaling one of the two extreme eigenvectors makes
		their nonzero pairing equal to $r$. Hence there exists
		$Q\in\GL_4(K_A)$ such that
		\[
		Q^*H_KQ=r\bigl(\mathbb H\perp\langle1,s/r\rangle\bigr).
		\]
		This gives an extraction with the same $A$ and with coefficient
		$q_A=s/r$. Taking $P=P_0Q$ and $a=br\in M_0^\times$, we obtain
		\[
		P^*H_\tau P=aH_A.
		\]
		The conjugated group is $k_A$-Zariski-dense in $\SU(H_A)$: a
		proper Zariski closure defined over $k_A$ would remain proper
		after extension to $M$, contradicting the density of the original
		representation after splitting. Taking determinants and using
		$\det H_\tau=-\tau$ and $\det H_A=-q_A$ gives
		\[
		-\tau\,N_{M/M_0}(\det P)=-a^4q_A,\qquad
		q_A=\tau\,N_{M/M_0}\!\left(\frac{\det P}{a^2}\right).
		\]
		
		The same calculation applies to every other normalized model of
		this representation. Indeed, two models over $M$ of the same
		complex representation admit an invertible complex intertwiner.
		The intertwining equations are linear over $M$, a finite
		subsystem suffices, and the solution space is one-dimensional
		by absolute irreducibility. This line has a generator over $M$,
		which is necessarily invertible. The intertwiner can therefore
		be chosen over $M$, and the invariant Hermitian forms differ by
		a scalar in $M_0^\times$. The identity applies to the coefficient
		of the given extraction, with no subsequent choice depending on
		an additional word.
	\end{proof}
	
	\subsection{The obstruction for every extraction and every finite family}
	
	\begin{proof}[Proof of Theorem~\ref{thm:transcendental-obstruction}]
		Lemma~\ref{lem:uniform-seed} and Proposition~\ref{lem:unitary-deformation}
		provide $B_1,\dots,B_4$, the finitely generated field
		$F\subset\mathbb R$, and a common interval $(0,\delta)$.
		The field $F$ and its real algebraic closure are countable, so
		there exists $T\in(0,\min\{1,\delta\})$ transcendental over $F$.
		Fix this parameter and the group
		\[
		\Lambda=\langle g_1(T),\ldots,g_4(T)\rangle
		\subset\SU(H_T),\qquad H_T=\mathbb H\perp\langle1,T\rangle.
		\]
		The proposition supplies all the stated geometric properties and
		regularity of every nontrivial word and power. The matrices have
		entries in $M_{\mathrm b}=F(i,T)$.
		
		Now fix any $A\in\Lambda\setminus\{1\}$ and put
		\[
		M_A=M_{\mathrm b}(E(A)),\qquad
		K_A=\mathbb Q(\operatorname{tr}\Lambda,E(A)),\qquad k_A=K_A^c.
		\]
		The conjugate-reciprocal characteristic polynomials make $M_A/M_{\mathrm b}$
		finite and $c$-stable, with $K_A\subset M_A$. By
		Lemma~\ref{lem:global-fields}, density persists from $F(T)$ to $M_A^c$.
		Since $M_A=M_A^c(i)$, this is base change of the same unitary form.
		
		Apply Proposition~\ref{prop:intrinsic-discriminant} with $M=M_A$ and
		$\tau=T$. It gives an extraction over the exact field $K_A$,
		with $A$ diagonal, and the coefficient of this extraction satisfies
		\begin{equation}\label{eq:extracted-norm-class}
			P_A^*H_TP_A=a_AH_A,\qquad
			q_A=T u_Ac(u_A),\qquad
			u_A=\frac{\det P_A}{a_A^2}\in M_A^\times.
		\end{equation}
		This identity holds for every normalized extraction by the
		intertwiner descent in that proposition. Finite generation of
		$K_A$ follows from Lemma~\ref{lem:trace-algebra}; the proposition
		also proves that $c$ is nontrivial and $[K_A:k_A]=2$.
		
		Let $\gamma_1,\ldots,\gamma_r$ be any finite family of nontrivial
		elements of $\Lambda$, allowing $r=0$. Apply
		Lemma~\ref{lem:conjugation-compatible-lifting} over $E'=E=F(i)$ to
		\[
		\chi_A(t,X),\quad\chi_{\gamma_1}(t,X),\ldots,
		\chi_{\gamma_r}(t,X).
		\]
		Here each symbol denotes the characteristic polynomial of the
		corresponding word in the rational generators $g_j(t)$. All their
		coefficients are regular at zero and on $(0,\delta)$: the
		generator matrices and their inverses introduce only powers of
		$1+t$ as denominators. The reciprocal identity follows from
		unitarity, separability at zero from \eqref{eq:triangular-reduction} and
		Lemma~\ref{lem:uniform-seed}, and the loxodromic root pattern on a
		common interval from Proposition~\ref{lem:unitary-deformation}.
		The spectrum of $A$ may vary with $t$; it is therefore included
		among the polynomials, rather than adjoined to the constant field.
		
		The final field in the lifting argument is
		\[
		\mathscr L_{A,\boldsymbol\gamma}
		=M_A(E(\gamma_1),\ldots,E(\gamma_r)).
		\]
		It admits an integer-valued valuation $w$, trivial on
		$F(i)^\times$, with $w(T)=1$ and $w\circ c=w$.
		Equation~\eqref{eq:extracted-norm-class} gives
		\begin{equation}\label{eq:odd-valuation}
			w(q_A)=1+2w(u_A)\in2\mathbb Z+1.
		\end{equation}
		If $q_A=zc(z)$ in
		$L_{A,\boldsymbol\gamma}=K_A(E(\gamma_1),\ldots,E(\gamma_r))$,
		the inclusion of this field in $\mathscr L_{A,\boldsymbol\gamma}$
		and agreement of their complex conjugations would imply
		$w(q_A)=2w(z)$, a contradiction. Thus $q_A$ is not a norm in
		any finite spectral compositum, including the case with no
		additional spectra.
		
		Every trace is a sum of characteristic roots, so
		$L_{A,\infty}=\Q(E(g):g\in\Lambda)$, independently of $A$.
		This equality concerns the total field, without identifying the fields $K_A$.
		It is the directed union of the $c$-stable fields
		$L_{A,\boldsymbol\gamma}$. Lemma~\ref{prop:finite-level-obstruction}
		excludes a norm in the union: any witness would already involve
		finitely many spectra. No single valuation on the union is required.
		
		Taking $r=0$, so that $A$ is the only word in the lifting lemma,
		gives such a valuation on $M_A$. A valuation trivial on
		$\mathbb Q^\times$ vanishes on every nonzero element algebraic
		over $\mathbb Q$. Indeed, if its value on such an element were
		positive, the constant term of its minimal polynomial would be
		the unique term of least value; if the value were negative, the
		leading term would be the unique term of least value. Both
		possibilities contradict the ultrametric inequality. Since the value in
		\eqref{eq:odd-valuation} is odd, $q_A$ is transcendental over
		$\mathbb Q$, and $K_A$ has positive transcendence degree.
		
		Finally, let $D_T=\operatorname{diag}(1,1,1,\sqrt T)$ and
		$J=\mathbb H\perp\langle1,1\rangle$. Since $D_T^*JD_T=H_T$,
		the group $\Gamma=D_T\Lambda D_T^{-1}$ lies in $\SU(J)$ and
		retains all the geometric properties established above.
		Conjugation preserves traces and spectra, hence the intrinsic fields
		and the total spectral field. Every extraction from $\Gamma$ gives
		the same complex representation, so the last assertion of
		Proposition~\ref{prop:intrinsic-discriminant} applies to its coefficient $q_A$.
		A constant real change of basis takes $J$ to the usual diagonal form
		of signature $(3,1)$, proving the assertions for this single group.
	\end{proof}
	
	\begin{corollary}
		\label{cor:no-total-unitary-descent}
		For the group in Theorem~\ref{thm:transcendental-obstruction} and every $A\ne1$,
		no change of basis over $\mathbb C$ realizes the same
		representation in
		\[
		\SU(J)(L_{A,\infty}^c),\qquad
		J=\mathbb H\perp\langle1,1\rangle.
		\]
		As in the notation of the introduction, $\SU(J)(L^c)$ consists of
		the determinant-one matrices over $L$ preserving $J$ under actual
		complex conjugation, relative to the extension $L/L^c$.
		This conclusion concerns the extension of the global assertion
		to fields of positive transcendence degree. The intrinsic field
		of this group is not imaginary quadratic, so the counterexample
		does not fall under that number-field hypothesis.
	\end{corollary}
	
	\begin{proof}
		Write $L=L_{A,\infty}$ and take the extracted representation
		$\rho$ over $K_A\subset L$, preserving $H_A$. If an equivalent
		representation $\rho'$ had matrices in $\SU(J)(L^c)$, the
		intertwining equations
		\[
		Q\rho(\lambda)=\rho'(\lambda)Q\qquad(\lambda\in\Lambda)
		\]
		would form a linear system over $L$. A finite subsystem suffices
		because there are sixteen unknowns. Its complex solution space
		is one-dimensional by absolute irreducibility and contains an
		invertible matrix. Linear descent gives a generator
		$Q\in\GL_4(L)$, necessarily invertible. Uniqueness of the
		invariant form, together with Hermitian symmetry, gives
		\[
		Q^*JQ=aH_A\quad\text{with }a\in(L^c)^\times.
		\]
		Taking determinants and using $\det J=-1$ and
		$\det H_A=-q_A$, we obtain
		\[
		N_{L/L^c}(\det Q)=a^4q_A,\qquad
		q_A=N_{L/L^c}\!\left(\frac{\det Q}{a^2}\right),
		\]
		contrary to \eqref{eq:total-spectral-nonnorm}. The initial complex change of
		basis was not assumed to be defined over $L$ or to have
		determinant one; descent and the determinant identity account
		for both possibilities.
	\end{proof}
	
	\section*{Disclosure statement, funding, declarations, and acknowledgments}
	
	AI-assisted tools were used for language editing, organization, and consistency checks. The authors are fully responsible for all mathematical statements, proofs, computations, citations, and the final wording of the article.
	
	This work was partially supported, for A. Cano, by SECIHTI--SNII grant 104023 and PAPIIT--UNAM project IN112424; and, for M. M. López-Flores, by the project EVAAT\_GCN, FINEP agreement 01.24.0661.00, ref. 3311/24. 
	
	All authors contributed equally to this work and declare no competing interests. The authors thank the staff of UCIM--UNAM and the Instituto Militar de Engenharia (IME--RJ, Brazil) for their hospitality during the preparation of this article. They are also grateful to John Parker for valuable conversations related to complex hyperbolic geometry and the themes of this work.
	
	\noindent\textbf{Code and data availability.}
	The computational supplement, verification programs, input data, and recorded
	outputs accompanying this article are associated with the Zenodo record
	\cite{SpectralFieldsComputationsZenodo}.

	% BEGIN AUTHOR INFORMATION
	\bigskip
	\begingroup
	\small
	\noindent\textsc{\'Angel Cano}\\
	Instituto de Matem\'aticas, Unidad Cuernavaca,\\
	Universidad Nacional Aut\'onoma de M\'exico,\\
	Cuernavaca, Morelos, Mexico\\
	\textit{E-mail address}: \href{mailto:angelcano@im.unam.mx}{\nolinkurl{angelcano@im.unam.mx}}
	
	\medskip
	\noindent\textsc{Nikolay Gusevskii}\\
	Departamento de Matem\'atica, Instituto de Ci\^encias Exatas,\\
	Universidade Federal de Minas Gerais,\\
	Belo Horizonte, MG, Brazil\\
	\textit{E-mail address}: \href{mailto:nikolay@mat.ufmg.br}{\nolinkurl{nikolay@mat.ufmg.br}}
	
	\medskip
	\noindent\textsc{Marlon M. L\'opez-Flores}\\
	Laborat\'orio de Intelig\^encia Artificial, Rob\'otica e Cibern\'etica (LIARC),\\
	Instituto Militar de Engenharia (IME-RJ),\\
	Rio de Janeiro, Brazil\\
	\textit{E-mail address}: \href{mailto:marlon.michael@ime.eb.br}{\nolinkurl{marlon.michael@ime.eb.br}}
	\par
	\endgroup
	% END AUTHOR INFORMATION
\end{document}